\documentclass{article}%
\usepackage{graphicx}
\usepackage{amsmath}
\usepackage{amsfonts}
\usepackage{amssymb}%
\newtheorem{theorem}{Theorem}

\newtheorem{conjecture}[theorem]{Conjecture}

\newtheorem{definition}[theorem]{Definition}

\newtheorem{lemma}[theorem]{Lemma}

\newtheorem{problem}[theorem]{Problem}
\newtheorem{proposition}[theorem]{Proposition}
\newtheorem{remark}[theorem]{Remark}

\newenvironment{proof}[1][Proof]{\textbf{#1.} }{\ \rule{0.5em}{0.5em}}
\begin{document}

\title{On the asymptotic scalar curvature ratio of complete Type I-like ancient
solutions to the Ricci flow on noncompact 3-manifolds}
\author{Bennett Chow\\University of California, San Diego
\and Peng\ Lu\\University of Oregon}
\date{}
\maketitle

\section{Introduction}

\label{Dedicated to Peter Li in honor of his 10 years as Editor-in-Chief of Communications in Analysis and Geometry.}%

Complete noncompact Riemannian manifolds with nonnegative sectional curvature
arise naturally in the Ricci flow when one takes the limits of dilations about
a singularity of a solution of the Ricci flow on a compact $3$-manifold
\cite{H-95a}. To analyze the singularities in the Ricci flow one needs to
understand these manifolds in depth. There are three invariants, asymptotic
scalar curvature ratio, asymptotic volume ratio and aperture, that have been
used to study the geometry of these manifolds at infinity.

Let $\left(  \mathcal{M}^{n},g\right)  $ be a complete noncompact Riemannian
manifold with nonnegative sectional curvature and let $O\in\mathcal{M}$ be
some point which we call the origin. The \textit{asymptotic scalar curvature
ratio} (ASCR) is defined by
\begin{equation}
\text{ASCR}\left(  \mathcal{M},g\right)  =\limsup_{d\left(  x,O\right)
\rightarrow\infty}R\left(  x\right)  \cdot d\left(  x,O\right)  ^{2},
\label{ASCR definition}%
\end{equation}
where $R$ is the scalar curvature. ASCR is a measure of how non-flat the
manifold is at infinity. Since the sectional curvature is nonnegative, there
is a positive constant $c$ depending only on $n$ such that%
\begin{align*}
c^{-1}\cdot\limsup_{d\left(  x,O\right)  \rightarrow\infty}\left\vert
Rm\left(  x\right)  \right\vert \cdot d\left(  x,O\right)  ^{2}  &
\leq\text{ASCR}\left(  \mathcal{M},g\right) \\
&  \leq c\cdot\limsup_{d\left(  x,O\right)  \rightarrow\infty}\left\vert
Rm\left(  x\right)  \right\vert \cdot d\left(  x,O\right)  ^{2}%
\end{align*}
where $Rm$ is the Riemann curvature tensor. In the literature sometimes
\[
\limsup_{d\left(  x,O\right)  \rightarrow\infty}|Rm\left(  x\right)  |\cdot
d\left(  x,O\right)  ^{2}%
\]
is used as the definition of ASCR. It is clear that ASCR is independent of the
choice of origins and is invariant under scaling.

ASCR has been used to study gap theorems. In particular, in \cite{ESS-89}
Eschenburg, Schroeder and Strake proved that if $\left(  \mathcal{M}%
^{2k+1},g\right)  $ is a complete noncompact odd-dimensional Riemannian
manifold with positive sectional curvature, then ASCR$\left(  \mathcal{M}%
,g\right)  >0$. These types of results are generally referred to as gap
theorems since they show the existence of a gap between flat $\mathbb{R}^{n}$
and metrics of positive curvature on $\mathbb{R}^{n}$. Gap theorems have been
proved by Greene-Wu \cite{GW-82}, Kasue and Sugahara \cite{KS-87}, and Drees
\cite{D-94}. For a survey of the history of gap theorems related to the notion
of ASCR, see \cite{G-97}. ASCR has also been used to study the structure of
manifolds at infinity \cite{PT-01}.

Below we give a few examples about how ASCR is used in singularity analysis in
the Ricci flow.

\begin{enumerate}
\item In \cite{H-95a} Hamilton showed that for a solution to the Ricci flow on
a compact 3-manifold forming a Type II singularity and satisfying an
injectivity radius estimate, there exists a sequence of dilations converging
to a complete solution $\left(  \mathcal{M}_{\infty}^{3},g_{\infty}\left(
t\right)  \right)  $ defined for all $t\in\left(  -\infty,\infty\right)  $
with nonnegative bounded sectional curvature and attaining its maximum of the
scalar curvature on space and time. In \cite{H-93b} it is shown using the
differential Harnack inequality of Li-Yau-Hamilton type \cite{H-93a} that the
universal covering solution $\left(  \widetilde{\mathcal{M}}_{\infty}%
^{3},\tilde{g}_{\infty}\left(  t\right)  \right)  $ of such a solution must be
a stationary solution (also called \textit{Ricci solitons}) of the Ricci flow
in the space of metrics modulo diffeomorphisms flowing along a gradient vector
field, that is, there exists a 1-parameter family of diffeomorphisms
$\varphi_{t}:\widetilde{\mathcal{M}}_{\infty}^{3}\rightarrow\widetilde
{\mathcal{M}}_{\infty}^{3}$ such that $\tilde{g}_{\infty}\left(  t\right)
=\varphi_{t}^{\ast}\left(  \tilde{g}_{\infty}\left(  0\right)  \right)  $ and
the 1-parameter family of vector fields $X\left(  t\right)  $ generated by
$\varphi_{t}$ are the gradients of functions $f\left(  t\right)  .$ In \S 20
of \cite{H-95a} it is shown that such a Ricci soliton must have ~ASCR$\left(
g\left(  t\right)  \right)  =\infty$ for all $t\in\left(  -\infty
,\infty\right)  $.

\item In \S 22 of \cite{H-95a} Hamilton showed using a geometric result about
bumps of curvature from \S 21 of \cite{H-95a} that if a complete solution to
the Ricci flow with bounded curvature has ASCR$\left(  g\left(  t\right)
\right)  =\infty$, then one can perform dimension reduction. In particular,
provided there is a local injectivity radius estimate, there exists a sequence
of points $\left\{  y_{\alpha}\right\}  _{\alpha\in\mathbb{N}}$ in
$\mathcal{M}^{3}$ with $\lim_{\alpha\rightarrow\infty}d_{g\left(  0\right)
}\left(  y_{\alpha},O\right)  =\infty$ such that the sequence of dilated
solutions on balls $\left\{  \left(  B_{g_{\alpha}\left(  0\right)  }\left(
y_{\alpha},r_{\alpha}\right)  ,g_{\alpha}\left(  t\right)  \right)
,y_{\alpha}\right\}  _{\alpha\in\mathbb{N}},$ where
\[
g_{\alpha}\left(  t\right)  =R\left(  g_{\alpha}\right)  \left(  y_{\alpha
},0\right)  g\left(  t\cdot R\left(  g_{\alpha}\right)  ^{-1}\left(
y_{\alpha},0\right)  \right)
\]
and $\lim_{\alpha\rightarrow\infty}r_{\alpha}^{2}R\left(  g_{\alpha}\right)
\left(  y_{\alpha},0\right)  =\infty,$ converges to complete solution to the
Ricci flow $\left(  \mathcal{P}_{\infty},k_{\infty}\left(  t\right)
,y_{\infty}\right)  $ with bounded nonnegative sectional curvature which
splits metrically as the product of $\mathbb{R}$ and a solution on a surface
with positive curvature.\textrm{ }This is why it is important to know when a
noncompact ancient solution has infinite asymptotic scalar curvature ratio,
when this is true one should be able to perform dimension reduction.

\item In \S 5.3 of \cite{H-97} it is shown that for an ancient solution
$\left(  \mathcal{M},g\left(  t\right)  \right)  $ to the Ricci flow which is
complete with bounded positive curvature operator and satisfies certain
pinching conditions, then ASCR$\left(  \mathcal{M},g\left(  t\right)  \right)
=\infty$ if and only if for some fixed time the solution metric has an
arbitrarily necklike end (see \S 2 below for a definition).
\end{enumerate}

The \textit{asymptotic volume ratio} is defined by
\[
~AVR\left(  g\right)  \doteqdot\lim_{r\rightarrow\infty}\frac{\text{Vol}%
\left[  B\left(  O,r\right)  \right]  }{r^{n}},
\]
where $B(O,r)=\{x:d(O,x)<r\}$. It is also used in singularity analysis in the
Ricci flow. In particular, in \S 19 of \cite{H-95a} Hamilton showed that a
complete ancient Type I-like solution to the Ricci flow with bounded positive
curvature operator and finite ASCR must have AVR$\left(  g\left(  t\right)
\right)  >0$ for all $t$. He also showed there that the scalar curvature
decays exactly quadratically:%
\[
0<c\leq R\left(  x,t\right)  d_{g\left(  t\right)  }\left(  x,O\right)
^{2}\leq C
\]
where $c$ and $C$ depend on $O$ and $t.$

This paper arises from exploring the following idea inspired by the work of
Hamilton \cite{H-97}. \emph{If a piece of a positively curved complete
noncompact manifold is sufficiently close to a long standard cylinder, then
its asymptotic scalar curvature ratio is large.}

\section{Main results}

To state the main result of this paper we need to define the so-called
$\left(  \varepsilon,k,L\right)  $\emph{-}necks. First we recall a few basic
definitions concerning necks from \S 3.2 of \cite{H-97}. A topological
\textit{neck} in a differentiable manifold $\mathcal{M}^{n}$ is a local
diffeomorphism%
\[
N:S^{n-1}\times\left[  a,b\right]  \rightarrow\mathcal{M}^{n},
\]
for some $a<b$.

Let $g$ be a Riemannian metric on $\mathcal{M}$ and $g\doteqdot N^{\ast
}\left(  g\right)  $ be the pulled-back metric on on $S^{n-1}\times\left[
a,b\right]  .$ A neck is called \textit{normal} if it satisfies the following
five conditions.

\begin{enumerate}
\item (\textbf{good slices}) Each slice $S^{n-1}\times\left\{  z\right\}
\subset S^{n-1}\times\left[  a,b\right]  ,$ $z\in\left[  a,b\right]  ,$ has
constant mean curvature with respect to $g$.

\item (\textbf{good parametrizations of the slices}) The identity map
\[
\iota:\left(  S^{n-1}\times\left\{  z\right\}  ,\bar{g}\right)  \rightarrow
\left(  S^{n-1}\times\left\{  z\right\}  ,g\right)
\]
is harmonic for all $z\in\left[  a,b\right]  ,$ where $\bar{g}$ is the
standard metric.

\item (\textbf{taking in account that conformal maps of }$S^{2}$\textbf{ are
harmonic}) When $n=3$ the center of mass of $S^{2}\times\left\{  z\right\}
\subset\mathbb{R}^{3}\times\left\{  z\right\}  $ with respect to $g$ is the
origin:%
\[
\int_{S^{2}\times\left\{  z\right\}  }\vec{x}~dA_{g}\left(  \vec{x}\right)
=\vec{0}\in\mathbb{R}^{3}%
\]
for all $z\in\left[  a,b\right]  $, where $dA_{g}$ is the area form.

\item (\textbf{good spacing of the slices}) The spacing of the slices are
normalized using volume by%
\[
\text{Vol}\left(  S^{n-1}\times\left[  z,w\right]  ,g\right)  =\omega
_{n-1}\int_{z}^{w}r\left(  y\right)  ^{n}dy,
\]
where%
\[
r\left(  y\right)  \doteqdot\left(  \frac{\text{Area}\left(  S^{n-1}%
\times\left\{  y\right\}  ,g\right)  }{\omega_{n-1}}\right)  ^{\frac{1}{n-1}}%
\]
is called the \textit{mean radius}.

\item (\textbf{aligning the parametrizations}) If $\bar{V}$ is a Killing
vector field on a slice $\left(  S^{n-1}\times\left\{  z\right\}  ,\bar
{g}\right)  $, then the unit normal vector field $\nu$ of $S^{2}\times\left\{
z\right\}  \subset\mathbb{R}^{3}\times\left\{  z\right\}  $ with respect to
$g$ satisfies%
\[
\int_{S^{n-1}\times\left\{  z\right\}  }\bar{g}(\bar{V},\nu)dA_{\bar{g}}=0
\]
for all $z\in\left[  a,b\right]  $.\textrm{ }
\end{enumerate}

Note that in\ condition 4 it is important that the power of $r(y)$ is $n,$ not
$n-1$.

\bigskip\noindent A neck is called $\left(  \varepsilon,k\right)
$\textit{-cylindrical} if

\begin{enumerate}
\item (\textbf{conformally close to cylindrical in }$C^{0}$) the metric%
\[
\widehat{g}\left(  \theta,z\right)  \doteqdot r\left(  z\right)  ^{-2}g\left(
\theta,z\right)
\]
satisfies%
\[
\left\vert \widehat{g}-\bar{g}\right\vert _{\bar{g}}\leq\varepsilon,
\]

\item (\textbf{conformally close to cylindrical up to }$C^{k}$)%
\[
\left\vert \bar{\nabla}^{j}\widehat{g}\right\vert _{\bar{g}}\leq\varepsilon
\]
for all $1\leq j\leq k,$ where $\bar{\nabla}$ is the covariant derivative with
respect to $\bar{g},$

\item (\textbf{mean radius changes slowly})%
\[
\left\vert \frac{d^{j}\log r\left(  z\right)  }{dz^{j}}\right\vert
\leq\varepsilon
\]
for all $1\leq j\leq k.$
\end{enumerate}

A neck $N$ in a manifold $\mathcal{M}^{n}$ is called an \textit{embedded neck}
if $N:S^{n-1}\times\left[  a,b\right]  \rightarrow\mathcal{M}^{n}$ is an embedding.

\begin{definition}
(i) A map $N:S^{n-1}\times\left[  a,b\right]  \rightarrow\mathcal{M}$ is
called an $\left(  \varepsilon,k,L\right)  $\emph{-neck} if it is a normal
$\left(  \varepsilon,k\right)  $-cylindrical neck and $b-a\geq2L.$

(ii) A complete Riemannian manifold $\left(  \mathcal{M},g\right)  $ with
exactly one topological end is called to have an \textit{arbitrarily necklike
end} if for every $(\varepsilon,k,L)$ there exists an $(\varepsilon,k,L)$-neck
in $\left(  \mathcal{M},g\right)  $.
\end{definition}

Definition (ii) is from \S 5.3 of \cite{H-97} (p. 71). We write $\left(
\varepsilon_{1},k_{1},L_{1}\right)  \preccurlyeq\left(  \varepsilon_{2}%
,k_{2},L_{2}\right)  $ if $\varepsilon_{1}\leq\varepsilon_{2},$ $k_{1}\geq
k_{2}$ and $L_{1}\geq L_{2}$. In the next lemma we collect a few simple
properties of $(\varepsilon,k,L)$-necks.

\begin{lemma}
(i) $\left(  \varepsilon,k,L\right)  $-necks are scale-invariant, i.e., if
$N:S^{n-1}\times\left[  a,b\right]  \rightarrow(M^{n},g)$ is an $\left(
\varepsilon,k,L\right)  $-neck then $N:S^{n-1}\times\left[  a,b\right]
\rightarrow(M^{n},\lambda^{2}\cdot g)$ is also an $\left(  \varepsilon
,k,L\right)  $-neck for any $\lambda>0$.

(ii) If $N:S^{n-1}\times\left[  a,b\right]  \rightarrow\mathcal{M}$ is an
$\left(  \varepsilon_{1},k_{1},L_{1}\right)  $-neck and $\left(
\varepsilon_{1},k_{1},L_{1}\right)  \preccurlyeq\left(  \varepsilon_{2}%
,k_{2},L_{2}\right)  $, then $N$ is an $\left(  \varepsilon_{2},k_{2}%
,L_{2}\right)  $-neck.

(iii) If $S_{r}^{n-1}$ is the sphere of radius $r,$ then for all
$\varepsilon>0$ and $k\in\mathbb{N}$ the product manifold $\mathcal{M}%
=S_{r}^{n-1}\times\lbrack-rL,rL]$ has an $\left(  \varepsilon,k,L\right)
$-neck given by the obvious map.
\end{lemma}

The main result of this paper is the following theorem which will be proved in
\S \ref{Main result}.

\begin{theorem}
\label{Main-ASCR} For every odd integer $n\geq3$ and every positive
$C_{0}<\infty$ there exists $\left(  \varepsilon_{0},k_{0},L_{0}\right)  $
such that if $\left(  \mathcal{M}^{n},g\right)  $ is a complete, orientable,
noncompact Riemannian manifold with bounded positive sectional curvature and
$\mathcal{M}$ contains an embedded $\left(  \varepsilon,k,L\right)  $-neck
with $\left(  \varepsilon,k,L\right)  \preccurlyeq\left(  \varepsilon
_{0},k_{0},L_{0}\right)  ,$ then%
\[
\text{ASCR}\left(  g\right)  \geq C_{0}.
\]

\end{theorem}

The main result may be thought of as a odd-dimensional \textit{quantitative}
version of the \textquotedblleft only if\textquotedblright\ part of Theorem
3.1 in \S 5.3 of \cite{H-97}, which states that a complete noncompact ancient
solution to the Ricci flow on a four-manifold with bounded positive curvature
operator satisfying certain pinching conditions has an arbitrarily necklike
end at sometime if and only if the asymptotic scalar curvature ratio is
infinite. In particular, we obtain a weaker characterization of ASCR$\left(
g\left(  t\right)  \right)  =\infty,$ which is one of the ingredients used in
the proof of Theorem \ref{Type I ancient have infinite ASCR} below.

The main result has the following useful consequence in dimension three.
Recall that an ancient solution is called \textit{Type-I like} if
\[
\sup_{\mathcal{M}^{n}\times(-\infty,0]}\left\vert t\right\vert \cdot\left\vert
Rm\left(  x,t\right)  \right\vert <\infty.
\]

\begin{theorem}
\label{Type I ancient have infinite ASCR}If $\left(  \mathcal{M}^{3},g\left(
t\right)  \right)  ,$ $-\infty<t<\omega,$ is a complete noncompact ancient
Type I-like solution to the Ricci flow with bounded positive sectional
curvature on an orientable $3$-manifold, then ~ASCR$\left(  g\left(  t\right)
\right)  =\infty$ for all $t\in\left(  -\infty,\omega\right)  .$
\end{theorem}

We shall prove this theorem in
\S \ref{Existence of necklike points in ancient solutions}. This answers a
conjecture of Hamilton when $n=3$. In \S 22 of \cite{H-95a} (p. 93) he writes:

\begin{quotation}
We do not know any examples of complete noncompact ancient solutions of
positive curvature operator with $Rs^{2}<\infty$ and $R|t|<\infty$, and we
conjecture none exist, since the curvature has had plenty of space and time to dissipate.
\end{quotation}

\begin{conjecture}
There does not exist complete noncompact ancient Type I-like solutions to the
Ricci flow with bounded positive curvature operator.
\end{conjecture}

By the theorem above such solutions, if they were to exist, must necessarily
have infinite asymptotic scalar curvature ratio. When $n=2$ it is proved in
Theorem 26.1 of \cite{H-95a} that there does not exist complete noncompact
ancient Type I-like solutions to the Ricci flow with positive sectional
curvature. For the mean curvature flow, which in general has very similar
properties as the Ricci flow, there are no complete noncompact strictly convex
Type I-like ancient solutions by Huisken's classification using his
monotonicity formula \cite{Hu-90}.

\section{Relative volumes and necks\label{rel vol and neck}}

The main result of this section, Proposition \ref{rel vol comp 1} below, is to
combine the relative volume comparison theorem (Lemma
\ref{zhu's version of rel vol comp thm}) and the existence of an embedded
$\left(  \varepsilon,k,L\right)  $-neck to show that there are small relative
volumes. Note that this result holds in both odd and even dimensions.

\subsection{Busemann functions\label{busemann functions}}

Let $(\mathcal{M}^{n},g)$ be a complete noncompact Riemannian manifold. Given
a point $Q\in\mathcal{M}$ and a ray $\gamma$ emanating from $Q$, the Busemann
function $b_{\gamma}:M\rightarrow\mathbb{R}$ associated to $\gamma$ is defined
by
\[
b_{\gamma}(x)\doteqdot\lim_{t\rightarrow+\infty}[t-d(\gamma(t),x)].
\]
Let $\mathcal{R}$ be the set of all rays emanating from $Q$. The Busemann
function $b_{Q}:M\rightarrow\mathbb{R}$ with the base point $Q$ is defined by
\[
b_{Q}(x)=\sup_{\gamma\in\mathcal{R}}b_{\gamma}(x)
\]
We collect some well-known properties of Busemann functions on complete
noncompact Riemannian manifolds with nonnegative sectional curvature in the
following lemma (for a proof see, for example, \cite{CG-72} and \cite{LT-87}).

\begin{lemma}
\label{lem prop busemann} Let $(\mathcal{M}^{n},g)$ be a complete noncompact
Riemannian manifold with nonnegative sectional curvature. Then

(i) The Busemann function $b_{Q}$ is proper, Lipschitz with Lipschitz constant
$1$ and bounded from below.

(ii) The sublevel sets $C_{r}\doteqdot\left\{  x\in\mathcal{M}:b_{Q}\left(
x\right)  \leq r\right\}  $ are compact and totally convex.

(iii) Let $S_{r}\doteqdot b_{Q}^{-1}(r)$ denote level set of the Busemann
function. Then for any $r_{1}\leq r_{2}$%
\[
S_{r_{1}}=\left\{  x\in C_{r_{2}}:d\left(  x,S_{r_{2}}\right)  =r_{2}%
-r_{1}\right\}  .
\]

(iv) $|b_{Q}\left(  x\right)  |\leq d\left(  x,Q\right)  $, which implies
$B\left(  Q,r\right)  \subset C_{r}.$
\end{lemma}

The following lemma says that any sufficiently long minimal geodesic segment
emanating from $Q$ can be well approximated by a ray emanating from $Q$ and
that the Busemann function is similar at infinity to the distance function to
the base point $Q$ (see Li-Tam \cite{LT-87} Theorem 2.3 on p. 177, Kasue
\cite{K-88} Lemma 1.4 on p. 598, Drees \cite{D-94} Lemma 1 on p. 80 and for an
exposition see \cite{CKE}). For the convenience of the reader we give a proof here.

\begin{lemma}
\label{Busemann fcn sublevel containments} Let $(\mathcal{M}^{n},g)$ be a
complete noncompact Riemannian manifold with nonnegative sectional curvature.
Define $\theta:[0,\infty)\rightarrow\lbrack0,\pi]$ by%
\[
\theta(r)=\sup_{\sigma\in\mathcal{S}(r)}\inf_{\gamma\in\mathcal{R}%
}\measuredangle_{Q}(\sigma^{\prime}(0),\gamma^{\prime}(0)),
\]
where $\mathcal{S}(r)$ is the set of all minimal geodesic segments $\sigma$ of
length $L(\sigma)\geq r$ emanating from $Q$. Then

(i) $\theta(r)$ is a nonincreasing function of $r$ and
\[
\lim_{r\rightarrow+\infty}\theta\left(  r\right)  =0.
\]

(ii) The function $b_{Q}$ and $d\left(  \cdot,Q\right)  $ are asymptotically
equal, more precisely%
\begin{equation}
\left(  1-\theta\left[  d\left(  x,Q\right)  \right]  \right)  \cdot d\left(
x,Q\right)  \leq b_{Q}\left(  x\right)  \leq d\left(  x,Q\right)
\label{Busemann fcn lower bound dist}%
\end{equation}
for all $x\in\mathcal{M}$.
\end{lemma}

\begin{proof}
(i) From the definition it is clear that $\theta(r)$ is a nonincreasing
function of $r.$ If $\lim_{r\rightarrow\infty}\theta\left(  r\right)  \neq0$,
there exists $\epsilon>0,$ a sequence of points $p_{i}\in\mathcal{M}$ with
$d(p_{i},Q)\nearrow+\infty,$ and minimal geodesic segments $\sigma_{i}$
joining $Q$ and $p_{i}$ parametrized by arc length such that
\[
\measuredangle_{Q}(\sigma_{i}^{\prime}(0),\gamma^{\prime}(0))\geq\epsilon
\]
for each $i$ and all rays $\gamma\in\mathcal{R}$. By the compactness of the
unit sphere in $T_{Q}\mathcal{M}$, there is a subsequence such that
$\lim_{j\rightarrow+\infty}\sigma_{i_{j}}^{\prime}(0)\doteqdot V_{\infty}$
exists. Let $\sigma_{\infty}:[0,+\infty)\rightarrow\mathcal{M}$ be the unique
geodesic with $\sigma_{\infty}(0)=Q$ and $\sigma_{\infty}^{\prime
}(0)=V_{\infty}.$ It is clear that $\sigma_{\infty}\in\mathcal{R}$. In
particular the condition
\[
\measuredangle_{Q}(\sigma_{i_{j}}^{\prime}(0),\sigma_{\infty}^{\prime}%
(0))\geq\epsilon
\]
is impossible for large enough $j.$ (i) is proved

(ii) For any $x\in\mathcal{M}$ let $\sigma$ be a minimal geodesic from $Q$ to
$x$ and $r=d(x,Q)$. Since $\mathcal{R}$ is a closed set, there is a ray
$\gamma\in\mathcal{R}$ such that $\measuredangle_{Q}(\sigma^{\prime}%
(0),\gamma^{\prime}(0))\leq\theta(r)$. Then it follows from sectional
curvature $K_{\mathcal{M}}\geq0$ and Toponogov's comparison theorem that
$d(x,\gamma(r))\leq\theta(r)\cdot r.$ Hence
\[
b_{Q}(x)\geq b_{\gamma}(x)\geq r-d(x,\gamma(r))\geq(1-\theta(r))\cdot r.
\]
The lemma is proved.
\end{proof}

\subsection{Necks in manifolds with positive sectional curvature}

We call a neck $N:S^{n-1}\times\left[  a,b\right]  \rightarrow\mathcal{M}$
\textit{absolute }$\left(  \varepsilon,k\right)  $\textit{-cylindrical} if it
satisfies the following two inequalities.

\begin{enumerate}
\item
\begin{equation}
\left\vert \frac{1}{r^{2}\left(  \frac{a+b}{2}\right)  }\cdot g-\bar
{g}\right\vert _{\bar{g}}\leq\varepsilon\text{ \ on }S^{n-1}\times\left[
a,b\right]  \label{absol neck ineq 1}%
\end{equation}
where $r\left(  \frac{a+b}{2}\right)  $ is the mean radius of $S^{n-1}%
\times\{\frac{a+b}{2}\}$.

\item
\begin{equation}
\left\vert \bar{\nabla}^{j}\left(  \frac{1}{r^{2}\left(  \frac{a+b}{2}\right)
}\cdot g\right)  \right\vert _{\bar{g}}\leq\varepsilon\text{ \ on }%
S^{n-1}\times\left[  a,b\right]  \label{absol neck ineq 2}%
\end{equation}
for all $1\leq j\leq k,$ where $\bar{\nabla}$ is the covariant derivative with
respect to $\bar{g}$.
\end{enumerate}

We call a topological neck an \textit{absolute }$\left(  \varepsilon
,k,L\right)  $\textit{-neck} if it is a normal absolute $\left(
\varepsilon,k\right)  $-cylindrical neck and $b-a\geq2L$. When $\left(
\varepsilon_{1},k_{1},L_{1}\right)  \preccurlyeq\left(  \varepsilon_{2}%
,k_{2},L_{2}\right)  ,$ an absolute $\left(  \varepsilon_{1},k_{1}%
,L_{1}\right)  $-neck is an absolute $\left(  \varepsilon_{2},k_{2}%
,L_{2}\right)  $-neck. The following result, whose proof is given in Appendix
A, holds for necks in arbitrary Riemannian manifolds.

\begin{lemma}
\label{lem absolute neck} Given $\left(  \varepsilon,k,L\right)  ,$ there
exists $\varepsilon^{\prime}=\varepsilon^{\prime}(\varepsilon,k,L)\leq
\varepsilon$ such that if $N:S^{n-1}\times\left[  a,b\right]  \rightarrow
\mathcal{M}$ is an $\left(  \varepsilon^{\prime},k,L\right)  $-neck, then $N$
is an absolute $\left(  \varepsilon,k,L\right)  $-neck.
\end{lemma}

Let $N:S^{n-1}\times\left[  -L,L\right]  \rightarrow\mathcal{M}^{n}$ be an
embedded neck in a complete noncompact Riemannian manifold $(\mathcal{M}%
^{n},g)$ with positive sectional curvature. By Gromoll and Meyer \cite{GM-69},
$\mathcal{M}^{n}$ is diffeomorphic to $\mathbb{R}^{n}$. Thus it follows from
the solution of the Schoenflies Conjecture in dimension $n\neq4$ \cite{M-59},
\cite{B-60} that the center sphere $N\left(  S^{n-1}\times\left\{  0\right\}
\right)  $ bounds a differentiable ball in $\mathcal{M}^{n}$ when $n\neq4$.
When $n=4$ Hamilton proved the following lemma (see Theorem 1.1 in \S 7 of
\cite{H-97}, the proof there works for all $n\geq2$).

\begin{lemma}
\label{lem Schoenflies} There exists $\left(  \varepsilon_{a},k_{a}%
,L_{a}\right)  $ having the following property. For any $\left(
\varepsilon,k,L\right)  \preccurlyeq\left(  \varepsilon_{a},k_{a}%
,L_{a}\right)  $ and any complete noncompact Riemannian manifold
$(\mathcal{M}^{n},g)$ with positive sectional curvature which has an embedded
$\left(  \varepsilon,k,L\right)  $-neck $N:S^{n-1}\times\left[  -L,L\right]
\rightarrow\mathcal{M}^{n}$, the center sphere $N\left(  S^{n-1}\times\left\{
0\right\}  \right)  $ bounds a differentiable ball in $\mathcal{M}^{n}$.
\end{lemma}

\subsection{Necks and relative volumes}

Let $\omega_{n-1}$ be the volume of the sphere $S^{n-1}$ of radius $1$, and
$\varepsilon_{b}\doteqdot\varepsilon^{\prime}(\frac{1}{10},1,L_{b})$ be a
function of $L_{b}$\ as in Lemma\ \ref{lem absolute neck}.

\begin{proposition}
\label{prop rel vol comp}\label{rel vol comp 1} For any $\delta>0$ there is
$L_{b}\geq\max(L_{a},16)$ having the following property. For any $\left(
\varepsilon,k,L\right)  \preccurlyeq\left(  \min(\varepsilon_{a}%
,\varepsilon_{b}),k_{a},L_{b}\right)  $, if there is an embedded
$(\varepsilon,k,L)$-neck $N$ in a complete noncompact Riemannian manifold
$\left(  \mathcal{M}^{n},g\right)  $ with positive sectional curvature
\[
N:S^{n-1}\times\left[  -L,L\right]  \rightarrow(\mathcal{M}^{n},g),
\]
then without any loss of generality we may assume that the component of
$\mathcal{M}-N\left(  S^{n-1}\times\left[  -L_{b},L_{b}\right]  \right)  $
bounded by $N(S^{n-1}\times\{-L_{b}\})$ is diffeomorphic to a ball and for any
$Q\in N(S^{n-1}\times\{-L_{b}\})$ and there exists $r_{0}>0$\textrm{ }such
that for any $R_{2}\geq R_{1}\geq r_{0}$ the following relative volume
estimate holds:%
\begin{equation}
\frac{\text{Vol}\left[  ~\overline{B}(Q,R_{2})\backslash B(Q,R_{1})\right]
}{\frac{\omega_{n-1}}{n}\cdot(R_{2}^{n}-R_{1}^{n})}\leq\delta
.\label{relative vol is less than delta}%
\end{equation}

\end{proposition}

\begin{remark}
The $r_{0}>0$\textrm{ }we shall choose has the property that for all
$R_{2}\geq R_{1}\geq r_{0},$ $\overline{B}(Q,R_{2})\backslash B(Q,R_{1})$ is
contained in the component (diffeomorphic to $S^{n-1}\times\mathbb{R}$) of
$\mathcal{M}-N\left(  S^{n-1}\times\left[  -L_{b},L_{b}\right]  \right)  $
bounded by $N(S^{n-1}\times\{L_{b}\}).$
\end{remark}

Recall that a geodesic is called normal if it is parametrized by arc length.
To prove this proposition we need the following form of the relative volume
comparison theorem (see Theorem 3.1, p. 226 in Zhu \cite{Z-97}). Let $\Gamma$
be any measurable subset of the unit sphere $S_{p}^{n-1}\subset T_{p}%
\mathcal{M}$. Given $r\leq R,$ let%
\[
A_{r,R}^{\Gamma}(p)\doteqdot\left\{  x\in\mathcal{M}\text{ : }%
\begin{array}
[c]{c}%
r\leq d(x,p)\leq R\text{ and there exists a }\\
\text{normal minimal geodesic }\gamma\text{ from }\\
\gamma\left(  0\right)  =p\text{ to }x\text{ satisfying }\gamma^{\prime}%
(0)\in\Gamma
\end{array}
\right\}  .
\]
Fix any point $p_{H}$ in the simply-connected space form of dimension $n$ and
of constant sectional curvature $H,$ let $A_{r,R}^{\Gamma}(p_{H})$ be the
corresponding set in the space form. Clearly $A_{r,R}^{\Gamma}(p)\subset
\overline{B}(p,R)\backslash B(p,r)$, and if $\Gamma=S_{p}^{n-1}$ then
$A_{r,R}^{\Gamma}(p)=\overline{B}(p,R)\backslash B(p,r)$.

\begin{lemma}
\label{zhu's version of rel vol comp thm}\emph{(Bishop-Gromov relative volume
comparison theorem)} Let $(\mathcal{M}^{n},g)$ be a complete Riemannian
manifold with $Ric_{\mathcal{M}}\geq(n-1)H.$ If $r\leq R,$ $s\leq S,$ $r\leq
s,R\leq S$ and $\Gamma$ is as above, then%
\[
\frac{\text{\emph{Vol}}(A_{s,S}^{\Gamma}(p))}{\text{\emph{Vol}}^{H}%
(A_{s,S}^{\Gamma}(p_{H}))}\leq\frac{\text{\emph{Vol}}(A_{r,R}^{\Gamma}%
(p))}{\text{\emph{Vol}}^{H}(A_{r,R}^{\Gamma}(p_{H}))},
\]
where \emph{Vol}$^{H}$ is the volume in the space form.
\end{lemma}

\noindent\textbf{Proof of Proposition \ref{rel vol comp 1}}.\textbf{ }By
multiplying the metric $g$ by a positive constant if necessary,\textrm{ }we
may assume that the center sphere $N(S^{n-1}\times\{0\})$ has mean radius
$r(0)=1$. Note that $N$ remains an $\left(  \varepsilon,k,L\right)  $-neck
after the scaling and the desired estimate
(\ref{relative vol is less than delta}) does not change after the scaling. Let
$N^{b}=N|_{S^{n-1}\times\left[  -L_{b},L_{b}\right]  }:S^{n-1}\times\left[
-L_{b},L_{b}\right]  \rightarrow\mathcal{M}$. It is clear that $N^{b}$ is a
$(\varepsilon,k,L_{b})$-neck.

It follows from Lemma \ref{lem Schoenflies} that $\mathcal{M}-N^{b}\left(
S^{n-1}\times\left[  -L_{b},L_{b}\right]  \right)  $ has two components
$\mathcal{U}_{1}$ and $\mathcal{U}_{2},$ where $\mathcal{U}_{1}$ is
diffeomorphic to an open ball $B^{n}$ and without loss of generality we may
assume bounds $N^{b}\left(  S^{n-1}\times\{-L_{b}\}\right)  ,$\ and
$\mathcal{U}_{2}$ is diffeomorphic to $S^{n-1}\times\mathbb{R}$ and \ bounds
$N^{b}\left(  S^{n-1}\times\{L_{b}\}\right)  $. Let $Q\in N^{b}(S^{n-1}%
\times\{-L_{b}\})$. Define for $R_{2}>R_{1}>0$
\[
\Gamma=\left\{
\begin{array}
[c]{c}%
\gamma^{\prime}(0)\in S_{Q}^{n-1}\ \text{:}%
\begin{array}
[c]{c}%
\text{there is }x\text{ such that }R_{1}\leq d(Q,x)\leq R_{2}\\
\text{ and there is exactly one normal minimal}\\
\text{ geodesic }\gamma\text{ from }\gamma(0)=Q\text{ to }x\text{.}%
\end{array}
\text{ }%
\end{array}
\right\}  .
\]
Then $A_{R_{1},R_{2}}^{\Gamma}(Q)$\textrm{\ }is a subset of $\overline
{B}(Q,R_{2})\backslash B(Q,R_{1})$ and\textrm{ }%
\begin{equation}
\text{Vol}\left[  A_{R_{1},R_{2}}^{\Gamma}(Q)\right]  =\text{Vol}\left[
~\overline{B}(Q,R_{2})\backslash B(Q,R_{1})\right]  \label{annulus volume}%
\end{equation}
since $\left[  ~\overline{B}(Q,R_{2})\backslash B(Q,R_{1})\right]  \backslash
A_{R_{1},R_{2}}^{\Gamma}(Q)$ is contained in the set of cut locus points of
$Q$, which has measure zero.

Let $H=0$. The corresponding space form is Euclidean space, and%
\[
\text{Vol}^{0}(A_{R_{1},R_{2}}^{\Gamma}(0))=\frac{m(\Gamma)}{n}\cdot(R_{2}%
^{n}-R_{1}^{n})
\]
where $m(\Gamma)$ is the measure of $\Gamma$ in $S_{Q}^{n-1}$. We will apply
Lemma \ref{zhu's version of rel vol comp thm} to (\ref{annulus volume}) to
prove (\ref{relative vol is less than delta}). Thus we need to find for
comparison another relative volume which is less or equal to $\delta$. We will
find this other relative volume by using the embedded neck.

Choose $r_{0}$ large enough (depending on the neck $N^{b}$ and the manifold
$\mathcal{M}$), so that $B\left(  Q,r_{0}\right)  \supset\mathcal{U}_{1}\cup
N^{b}\left(  S^{n-1}\times\left[  -L_{b},L_{b}\right]  ]\right)  $. This
implies that if $R_{2}>R_{1}>r_{0},$ then $A_{R_{1},R_{2}}^{\Gamma}%
(Q)\subset\mathcal{U}_{2}$.

For any normal minimal geodesic $\gamma_{0}:\left[  0,\ell_{0}\right]
\rightarrow\mathcal{M}$ with $\gamma_{0}(0)=Q,$ $\gamma_{0}\left(  \ell
_{0}\right)  \in A_{R_{1},R_{2}}^{\Gamma}(Q)$ and $\gamma_{0}^{\prime}%
(0)\in\Gamma,$ then $\gamma_{0}$ will intersect $N^{b}\left(  S^{n-1}%
\times\{0\}\right)  $ at some (exactly one) point, say $\gamma_{0}(w_{0}),$
and we claim that
\begin{equation}
w_{0}\geq\frac{9}{10}L_{b}~. \label{claim 1}%
\end{equation}
To see the claim, since $\varepsilon_{b}\doteqdot\varepsilon^{\prime}(\frac
{1}{10},1,L_{b})$, by Lemma \ref{lem absolute neck} $N^{b}$ is an absolute
$\left(  \frac{1}{10},1,L_{b}\right)  $-neck and hence (here we use the
assumption that mean radius $r(0)=1$)
\begin{equation}
\frac{9}{10}\left\vert ~\bullet~\right\vert _{N_{\ast}^{b}\overline{g}}%
\leq\left\vert ~\bullet~\right\vert _{g}\leq\frac{11}{10}\left\vert
~\bullet~\right\vert _{N_{\ast}^{b}\overline{g}} \label{metric equiv relat}%
\end{equation}
where $N_{\ast}^{b}\overline{g}$ is the push-forward metric. From $\gamma
_{0}\left(  w_{0}\right)  \in N^{b}\left(  S^{n-1}\times\{0\}\right)  $ and
$Q=\gamma_{0}\left(  0\right)  \in N^{b}\left(  S^{n-1}\times\{-L_{b}%
\}\right)  $ we have
\[
d_{N_{\ast}^{b}\bar{g}}\left(  \gamma_{0}\left(  0\right)  ,\gamma_{0}\left(
w_{0}\right)  \right)  \geq L_{b}.
\]
Hence
\[
w_{0}=d_{g}\left(  \gamma_{0}\left(  0\right)  ,\gamma_{0}\left(
w_{0}\right)  \right)  \geq\frac{9}{10}d_{N_{\ast}\bar{g}}\left(  \gamma
_{0}\left(  0\right)  ,\gamma_{0}\left(  w_{0}\right)  \right)  =\frac{9}%
{10}L_{b}.
\]

We now claim that (this is only a rough estimate) \textrm{ }%
\begin{equation}
w_{0}+10\leq r_{0}. \label{claim 2}%
\end{equation}
Since $\mathcal{M}\backslash B(Q,r_{0})\subset\mathcal{U}_{2},$ we have
\[
r_{0}\geq d_{g}\left(  Q,\mathcal{U}_{2}\right)  \geq d_{g}(Q,N^{b}%
(S^{n-1}\times\{L_{b}\}).
\]
Since any minimal geodesic from $Q$ to any $p\in N^{b}(S^{n-1}\times
\{L_{b}\})$ must intersect $N^{b}(S^{n-1}\times\{0\}),$ we have\textrm{ }%
\begin{align*}
d_{g}(Q,N^{b}(S^{n-1}\times\{L_{b}\})  &  \geq d_{g}(Q,N^{b}(S^{n-1}%
\times\{0\})\\
&  +d_{g}(N^{b}(S^{n-1}\times\{0\},N^{b}(S^{n-1}\times\{L_{b}\}).
\end{align*}
By (\ref{metric equiv relat}) the diameter
\[
diam_{g}(S^{n-1}\times\{0\})\leq\frac{11}{10}\cdot diam_{N_{\ast}^{b}%
\overline{g}}(S^{n-1}\times\{0\})=\frac{11\pi}{10},
\]
it follows%
\[
d_{g}(Q,N^{b}(S^{n-1}\times\{0\})\geq d_{g}(Q,\gamma(w_{0}))-\frac{11\pi}%
{10}=w_{0}-\frac{11\pi}{10}.
\]
>From (\ref{metric equiv relat}) we have
\begin{align*}
&  d_{g}(N^{b}(S^{n-1}\times\{0\},N^{b}(S^{n-1}\times\{L_{b}\})\\
&  \geq\frac{9}{10}\cdot d_{N_{\ast}^{b}\overline{g}}(N^{b}(S^{n-1}%
\times\{0\},N^{b}(S^{n-1}\times\{L_{b}\})=\frac{9}{10}\cdot L_{b}%
\end{align*}
\textrm{ }Hence%
\[
w_{0}-\frac{11\pi}{10}+\frac{9}{10}\cdot L_{b}\leq r_{0}.
\]
The claim follows since $L_{b}\geq16$.\textrm{ }

Now we choose the $r,R,s$ and $S$ in Lemma
\ref{zhu's version of rel vol comp thm} as
\[
r\doteqdot w_{0}-2,\;\;\;R\doteqdot w_{0}+2,\;\;\;s\doteqdot R_{1}%
,\;\;\;S\doteqdot R_{2}.
\]
It is clear from the choice of $L_{b}\geq16$ and $r_{0}<R_{1}$ that $r\leq
R\leq s\leq S$. We claim\bigskip

\noindent\textbf{Sublemma}%
\[
A_{w_{0}-2,w_{0}+2}^{\Gamma}(Q)\subset\left\{  x\in\mathcal{M}:d_{g}\left(
x,N^{b}\left(  S^{n-1}\times\{0\}\right)  \right)  \leq6\right\}  ,
\]
where $\Gamma$ is the set defined above.\smallskip

\begin{remark}
Intuitively, the set $\{x\in\mathcal{M}$ : $d_{g}\left(  x,N^{b}\left(
S^{n-1}\times\{0\}\right)  \right)  \leq6\}$ is close to a standard cylinder
of length $12$ and radius $1.$
\end{remark}

\begin{problem}
Given $R_{1}$ and $R_{2},$ let $\Gamma$ be defined as above.
\end{problem}

\begin{proof}
[Proof of sublemma]\label{added sentence 1}Let $R_{1},R_{2},$ $\Gamma$ and
$w_{0}$ be as above. Note that $\gamma_{0}(w_{0})\in N^{b}\left(
S^{n-1}\times\{0\}\right)  $ and $\Gamma$ corresponds to $R_{1}$ and $R_{2}.$
Given any point $x\in A_{w_{0}-2,w_{0}+2}^{\Gamma}(Q),$ let $\gamma:\left[
0,\ell_{1}\right]  \rightarrow\mathcal{M}$ be a normal minimal geodesic with
$\gamma(0)=Q,$ $\gamma^{\prime}(0)\in\Gamma,$ $\gamma\left(  \ell_{1}\right)
\in\overline{B}(Q,R_{2})\backslash B(Q,R_{1})$ and $\gamma\left(  \ell
_{x}\right)  =x$ for some $\ell_{x}\in\left[  w_{0}-2,w_{0}+2\right]  .$
\label{added sentence 2}The geodesic $\gamma$ exists since $x\in
A_{w_{0}-2,w_{0}+2}^{\Gamma}(Q)$ implies there exists a normal minimal
geodesic $\bar{\gamma}:\left[  0,\ell_{x}\right]  \rightarrow\mathcal{M}$ with
$\bar{\gamma}^{\prime}\left(  0\right)  \in\Gamma,$ $\bar{\gamma}\left(
\ell_{x}\right)  =x$ and $\ell_{x}\in\left[  w_{0}-2,w_{0}+2\right]  .$ Since
$\bar{\gamma}^{\prime}\left(  0\right)  \in\Gamma,$ $\bar{\gamma}$ extends to
a normal minimal geodesic $\gamma$ as above. From $\gamma\left(  \ell
_{1}\right)  \notin B(Q,R_{1}),$ we have $\ell_{1}\geq R_{1}>r_{0}.$ Then
$\gamma\left(  \ell_{1}\right)  \notin\mathcal{U}_{1}\cup N^{b}\left(
S^{n-1}\times\lbrack-L_{b},L_{b}]\right)  $ and hence $\gamma$ will intersect
$N^{b}\left(  S^{n-1}\times\{0\}\right)  $ at some point $\gamma\left(
w_{x}\right)  $ (for the same reason as before as applied to $\gamma_{0}$).
That is,
\[
\gamma(w_{x})\in N^{b}\left(  S^{n-1}\times\{0\}\right)  .
\]
Since the mean radius of\ $N^{b}\left(  S^{n-1}\times\{0\}\right)  $ is $1$,
we have by (\ref{metric equiv relat})
\[
d_{g}(\gamma(w_{x}),\gamma_{0}(w_{0}))\leq\frac{11\pi}{10}.
\]
Hence we get
\begin{equation}
\left\vert w_{x}-w_{0}\right\vert \leq\frac{11\pi}{10}\leq4. \label{Lx and L0}%
\end{equation}
>From the triangle inequality
\[
\left\vert d_{g}(Q,\gamma(w_{x}))-d_{g}(Q,\gamma_{0}(w_{0}))\right\vert \leq
d_{g}(\gamma(w_{x}),\gamma_{0}(w_{0})).
\]
On the other hand%
\begin{equation}
w_{0}-2\leq d_{g}(Q,x)=\ell_{x}\leq w_{0}+2. \label{L0 and lx}%
\end{equation}
Combining (\ref{Lx and L0}) and (\ref{L0 and lx}) we get
\begin{equation}
\left\vert l_{x}-w_{x}\right\vert \leq6, \label{Lx and lx}%
\end{equation}
which implies
\[
d_{g}(\gamma(w_{x}),x)=d_{g}\left(  \gamma\left(  w_{x}\right)  ,\gamma\left(
\ell_{x}\right)  \right)  \leq6.
\]
Since $\gamma(w_{x})\in N\left(  S^{n-1}\times\{0\}\right)  ,$
\[
d_{g}(N^{b}\left(  S^{n-1}\times\{0\}\right)  ,x)\leq6.
\]
This completes the proof of the sublemma.\medskip
\end{proof}

It follows from the sublemma and (\ref{metric equiv relat}) that
\begin{align*}
&  \text{Vol}\left[  A_{w_{0}-2,w_{0}+2}^{\Gamma}(Q)\right] \\
&  \leq\text{Vol}\left(  \{x\in\mathcal{M}:d_{g}(x,N^{b}\left(  S^{n-1}%
\times\{0\}\right)  )\leq6\}\right) \\
&  \leq\text{Vol}\left(  \{x\in\mathcal{M}:d_{N_{\ast}^{b}\overline{g}%
}(x,N^{b}\left(  S^{n-1}\times\{0\}\right)  )\leq\frac{10}{9}\cdot6\}\right)
\\
&  \leq\left(  \frac{11}{10}\right)  ^{n}\cdot\text{Vol}_{N_{\ast}%
^{b}\overline{g}}\left(  \{x\in\mathcal{M}:d_{N_{\ast}^{b}\overline{g}%
}(x,N^{b}\left(  S^{n-1}\times\{0\}\right)  )\leq\frac{10}{9}\cdot6\}\right)
\\
&  =\left(  \frac{11}{10}\right)  ^{n}\cdot\omega_{n-1}\cdot1^{n}\cdot
\frac{10}{9}\cdot12.
\end{align*}

We now finish the proof of the proposition. Applying Lemma
\ref{zhu's version of rel vol comp thm} we get
\[
\frac{\text{Vol}\left[  B(Q,R_{2})\backslash B(Q,R_{1})\right]  }%
{\frac{m(\Gamma)}{n}\cdot\left(  R_{2}^{n}-R_{1}^{n}\right)  }\leq
\frac{\text{Vol}\left(  A_{w_{0}-2,w_{0}+2}^{\Gamma}(Q)\right)  }%
{\frac{m(\Gamma)}{n}\cdot((w_{0}+2)^{n}-(w_{0}-2)^{n})}%
\]
since $w_{0}-2\leq w_{0}+2\leq r_{0}\leq R_{1}\leq R_{2}$. Replacing the
common factor $m(\Gamma)$ by $\omega_{n-1}$ and applying the estimate above,
we get%
\begin{align*}
\frac{\text{Vol}\left[  B(Q,R_{2})\backslash B(Q,R_{1})\right]  }{\frac
{\omega_{n-1}}{n}\cdot\left(  R_{2}^{n}-R_{1}^{n}\right)  }  &  \leq
\frac{\text{Vol}\left(  A_{w_{0}-2,w_{0}+2}^{\Gamma}(Q)\right)  }{\frac
{\omega_{n-1}}{n}\cdot((w_{0}+2)^{n}-(w_{0}-2)^{n})}\\
&  \leq\frac{\left(  \frac{11}{10}\right)  ^{n}\cdot\omega_{n-1}\cdot\frac
{10}{9}\cdot12}{\frac{\omega_{n-1}}{n}\cdot((w_{0}+2)^{n}-(w_{0}-2)^{n})}.
\end{align*}

We have proved by using (\ref{claim 1})
\begin{equation}
\frac{\text{Vol}\left[  B(Q,R_{2})\backslash B(Q,R_{1})\right]  }{\frac
{\omega_{n-1}}{n}\cdot\left(  R_{2}^{n}-R_{1}^{n}\right)  }\leq\frac{\left(
\frac{11}{10}\right)  ^{n}\cdot\omega_{n-1}\cdot\frac{10}{9}\cdot12}%
{\frac{\omega_{n-1}}{n}\cdot((\frac{9}{10}L_{b}+2)^{n}-(\frac{9}{10}%
L_{b}-2)^{n})}. \label{eq rel vol up est}%
\end{equation}
If we choose $L_{b}\geq\max(L_{a},16)$ satisfying
\begin{equation}
\frac{\left(  \frac{11}{10}\right)  ^{n}\cdot\omega_{n-1}\cdot\frac{10}%
{9}\cdot12}{\frac{\omega_{n-1}}{n}\cdot((\frac{9}{10}L_{b}+2)^{n}-(\frac
{9}{10}L_{b}-2)^{n})}\leq\delta, \label{eq for delta choice}%
\end{equation}
and choose $\varepsilon_{b}=\varepsilon^{\prime}(\frac{1}{10},1,L_{b})$, then
the proposition follows from (\ref{eq rel vol up est}) and
(\ref{eq for delta choice}).

\section{Proof of the main result \label{Main result}}

The main \ part of this section is devoted to estimate the relative volume in
(\ref{relative vol is less than delta}) from below by ASCR when $R_{1}$ and
$R_{2}$\ is large and dimension $n$ is odd; see Proposition
\ref{prop rel vol est below} below. The main result Theorem \ref{Main-ASCR} is
proved at the very end of this section.

\subsection{Asymptotic scalar curvature ratio\label{subsec ascr a,k func}}

Let $\left(  \mathcal{M}^{n},g\right)  $ be a complete, noncompact Riemannian
manifold with positive sectional curvature and $Q\in\mathcal{M}$. Define a
function $a:\mathbb{\bar{R}}^{+}\rightarrow\mathbb{\bar{R}}^{+}$ by
\begin{equation}
a(r)^{2}=\sup_{x\in\mathcal{M}\backslash B\left(  Q,r\right)  }R\left(
x\right)  d\left(  x,Q\right)  ^{2}, \label{eq for a(r)}%
\end{equation}
a function $\kappa(r):\mathbb{\bar{R}}^{+}\rightarrow\mathbb{\bar{R}}^{+}$ by%
\[
\kappa(r)=\sup_{x\in\mathcal{M}\backslash B\left(  Q,r\right)  }R\left(
x\right)  ,
\]
and a function $\rho(r):\mathbb{\bar{R}}^{+}\rightarrow\mathbb{\bar{R}}^{+}$
by
\[
\rho(r)\doteqdot\frac{\pi r}{4a(r)}.
\]

The following lemma is clear.

\begin{lemma}
\label{lem on a(r),k(r)} (i) $a(r)$ is positive and monotone nonincreasing.

(ii) ASCR$\left(  g\right)  =\lim_{r\rightarrow\infty}a(r)^{2}.$

(iii) $\kappa(r)\cdot r^{2}\leq a(r)^{2}.$
\end{lemma}

To prove Proposition \ref{prop rel vol est below}, we assume $ASCR(g)<+\infty$
since otherwise the proposition is clearly true. For any $\eta_{1}\in(0,1)$
there is $r_{1}=r_{1}(\eta_{1},\mathcal{M})$ such that \
\[
a(r)\leq\sqrt{ASCR(g)}+\eta_{1},\text{ for all }r\geq r_{1}.
\]

Clearly we have
\begin{align}
\frac{a(r)}{r}  &  \rightarrow0^{+}\;\;\;\text{as }r\rightarrow+\infty
,\nonumber\\
\rho(r)  &  \rightarrow+\infty\text{\ \ \ as }r\rightarrow+\infty\nonumber\\
\frac{r}{\rho(r)}  &  \leq\frac{4(\sqrt{ASCR(g)}+\eta_{1})}{\pi}%
\text{\ \ \ for }r\geq r_{1}. \label{r over rho}%
\end{align}

\subsection{The hypersurfaces $\widehat{S}_{r}(\rho)$%
\label{subsec hypersurface}}

Let $S_{r}$ be the level set of the Busemann function $b_{Q}$ (defined in
\S \ref{busemann functions}) and $C_{r}$ be the sublevel set of the Busemann
function $b_{Q}$.\ If $S_{r}$ is smooth, we define $\widehat{b}_{Q}\doteqdot
b_{Q}$, $\widehat{C}_{r}\doteqdot C_{r}$ and $\widehat{S}_{r}\doteqdot S_{r}.$
If $S_{r}$ is not smooth, since $\mathcal{M}$ has positive sectional curvature
so $b_{Q}$ is strictly convex \cite{GW-74}, we can smooth $b_{Q}$ (see, for
example, p. 158 of \cite{ESS-89} or \cite{GW-76}). For any positive $\eta
_{2}\leq1$ there is a smooth and strictly convex function $\widehat{b}_{Q}$
such that
\begin{equation}
\left\vert \widehat{b}_{Q}(x)-b_{Q}(x)\right\vert <\eta_{2}\text{ for all
}x\in\mathcal{M}. \label{eq for approx busemann}%
\end{equation}
We define $\widehat{C}_{r}\doteqdot\widehat{b}_{Q}^{-1}(-\infty,r)$ and
$\widehat{S}_{r}\doteqdot\widehat{b}_{Q}^{-1}(r)$. So in any case we have a
smooth and strictly convex hypersurface $\widehat{S}_{r}$ and a strictly
convex set $\widehat{C}_{r}$. It is clear that $\widehat{S}_{r}\doteqdot
\partial\widehat{C}_{r}$.

We define the hypersurfaces $\widehat{S}_{r}\left(  \rho\right)  $ parallel to
$\widehat{S}_{r}$ by%
\[
\widehat{S}_{r}\left(  \rho\right)  \doteqdot\left\{  x\in\mathcal{M}:d\left(
x,\widehat{C}_{r}\right)  =\rho\right\}  .
\]
The following, which is Lemma 2 on p. 157 of \cite{ESS-89}, gives an estimate
for the second fundamental form of these parallel hypersurfaces.

\begin{lemma}
\label{S_r(rho) is embedded etc}\emph{(bounds for the 2nd fundamental form of
parallel hypersurfaces)} (i) If sectional curvature $K\leq\varepsilon^{2}$ on
$\mathcal{M}-\widehat{C}_{r},$ then the parallel hypersurfaces are smooth
embedded hypersurfaces for $0<\rho<\pi/\left(  2\varepsilon\right)  $.

(ii) Let \thinspace$g_{r}^{\rho}\,$ and \thinspace$h_{r}^{\rho}\in C^{\infty
}\left(  S^{2}T^{\ast}\widehat{S}_{r}\left(  \rho\right)  \right)  $
\thinspace denote the first and second fundamental forms of $\widehat{S}%
_{r}\left(  \rho\right)  ,$ respectively. Then%
\[
-\varepsilon\tan\left(  \varepsilon\rho\right)  g_{r}^{\rho}\leq h_{r}^{\rho
}\leq\frac{1}{\rho}g_{r}^{\rho}.
\]

(iii) Taking $\rho=\pi/(4\varepsilon)$ in (ii) we have
\[
-\varepsilon g_{r}^{\pi/\left(  4\varepsilon\right)  }\leq h_{r}^{\pi/\left(
4\varepsilon\right)  }\leq\frac{4\varepsilon}{\pi}g_{r}^{\pi/\left(
4\varepsilon\right)  }.
\]

(iv) The Weingarten map $L_{r}^{\pi/\left(  4\varepsilon\right)  }%
:T\widehat{S}_{r}\left(  \pi/(4\varepsilon)\right)  \rightarrow T\widehat
{S}_{r}\left(  \pi/(4\varepsilon)\right)  $ satisfies%
\[
\left\Vert L_{r}^{\pi/\left(  4\varepsilon\right)  }\right\Vert \leq
\frac{4\varepsilon}{\pi},
\]
where $\left\Vert L\right\Vert \doteqdot\max_{\left\vert v\right\vert
=1}\left\vert L\left(  v\right)  \right\vert $.
\end{lemma}

\begin{remark}
Note that if $\varepsilon$ is small, then $\pi/\left(  2\varepsilon\right)  $
is large. That is, the parallel hypersurfaces $\widehat{S}_{r}\left(
\rho\right)  $ are smooth for large $\rho.$ Our conclusion is that, assuming
$K\leq\varepsilon^{2}$ on $\mathcal{M}-\widehat{C}_{r},$ the second
fundamental form of $\widehat{S}_{r}\left(  \pi/\left(  4\varepsilon\right)
\right)  $ is small. The reason we have the weaker $1/\rho$ upper bound on the
second fundamental form is that $\widehat{S}_{r}$ could be close to a point
(like a small sphere). In particular, the $1/\rho$ upper bound is sharp for
$\widehat{S}_{0}$ a point in euclidean space.
\end{remark}

To apply Lemma \ref{S_r(rho) is embedded etc} we need to estimate
$\varepsilon$ in terms of $\kappa(r)$ and hence $a(r)$ in regards to the
condition that the sectional curvature $K\leq\varepsilon^{2}$ on
$\mathcal{M}-\widehat{C}_{r}$ $.$ We will use the following elementary result
whose proof is given in Appendix B.

\begin{lemma}
\label{lem smoothed set} (i) $\widehat{S}_{r}\subset b_{Q}^{-1}((r-\eta
_{2},r+\eta_{2}))$ and $\widehat{S}_{r}\subset\mathcal{M}\backslash
\overline{B}(Q,r-\eta_{2}).$

(ii) $B(Q,r-\eta_{2})\subset b_{Q}^{-1}((-\infty,r-\eta_{2}))\subset
\widehat{C}_{r}\subset b_{Q}^{-1}((-\infty,r+\eta_{2})).$

(iii) $\widehat{S}_{r}(\rho)\subset\mathcal{M}\backslash\overline{B}%
(Q,r+\rho-\eta_{2})$.

(iv) If $\widehat{S}_{r}\subset B(Q,\eta)$, then $\widehat{S}_{r}(\rho)\subset
B(Q,\eta+\rho).$
\end{lemma}

\subsection{Estimate of the area of $\widehat{S}_{r}(\rho)$ by Gauss-Bonnet
formula}

Note that sectional curvatures are less or equal to scalar curvatures
pointwise on a manifold with nonnegative sectional curvature, i.e., $K_{x}\leq
R(x)$. By Lemma \ref{lem smoothed set} (ii) and Lemma \ref{lem on a(r),k(r)}%
\ (iii)
\[
\sup_{\mathcal{M}\backslash\widehat{C}_{r}}K_{x}~\leq\sup_{\mathcal{M}%
\backslash B\left(  Q,r-\eta_{2}\right)  }R(x)=\kappa(r-\eta_{2})\leq
\frac{a(r-\eta_{2})^{2}}{(r-\eta_{2})^{2}}=\left(  \frac{\pi}{4\rho\left(
r-\eta_{2}\right)  }\right)  ^{2}.
\]
Hence by Lemma \ref{S_r(rho) is embedded etc} (i), (ii) the hypersurface
$\widehat{S}_{r}(\rho)$ is smooth for $\rho\leq\rho(r-\eta_{2})$ and its
Weingarten map $L_{r}(\rho)$ is bounded by
\begin{equation}
\left\Vert L_{r}(\rho)\right\Vert \leq\frac{1}{\rho}\text{ for any }\rho
\leq\rho(r-\eta_{2}).\label{estimate for L_r(rho)}%
\end{equation}
Indeed, since $\tan\theta\leq\frac{1}{\theta}$ for $0<\theta\leq\pi/4,$ we
have $\frac{\pi}{4\rho\left(  r-\eta_{2}\right)  }\tan\left(  \frac{\pi}%
{4\rho\left(  r-\eta_{2}\right)  }\rho\right)  \leq\frac{1}{\rho}$ if
$\rho\leq\rho\left(  r-\eta_{2}\right)  $.

For the remainder of this section we assume that $n$ is odd and consider only
$\rho\in(0,\rho\left(  r-\eta_{2}\right)  )$. Let $m\doteqdot\left(
n-1\right)  /2$. We shall apply the Gauss-Bonnet formula to the hypersurface
$\widehat{S}_{r}(\rho)$ for $r\geq r_{1}$ as defined in subsection
\ref{subsec ascr a,k func}. Following \cite{ESS-89} \cite{GW-82}, let
$G_{r}\left(  \rho\right)  $ be the Gauss-Bonnet integrand of $\widehat{S}%
_{r}\left(  \rho\right)  $ with the induced metric. There are many instances
of the formula for $G_{r}\left(  \rho\right)  $ in the literature (see for
example p. 749 of \cite{C-44} or p. 740 of \cite{GW-82}), in general
$G_{r}\left(  \rho\right)  $ is defined by%
\[
\frac{2}{\omega_{n-1}}G_{r}\left(  \rho\right)  dV=\frac{1}{2^{2m}\pi^{m}%
m!}\sum\epsilon_{i_{1},\ldots,i_{2m}}\Omega_{i_{1}i_{2}}\wedge\cdots
\wedge\Omega_{i_{2m-1}i_{2m}},
\]
where $\Omega_{ij}$ are the curvature 2-forms of the induced metric on
$\widehat{S}_{r}\left(  \rho\right)  $. Recall that $\omega_{n-1}=$Vol$\left(
S^{2m}\right)  =2^{2m+1}\pi^{m}m!/\left(  2m\right)  !.$Hence%
\[
G_{r}\left(  \rho\right)  dV=\frac{1}{\left(  2m\right)  !}\sum\epsilon
_{i_{1},\ldots,i_{2m}}\Omega_{i_{1}i_{2}}\wedge\cdots\wedge\Omega
_{i_{2m-1}i_{2m}}.
\]
By the Gauss-Bonnet formula%
\begin{equation}
\frac{2}{\omega_{n-1}}\int_{\widehat{S}_{r}\left(  \rho\right)  }G_{r}\left(
\rho\right)  dV_{\widehat{S}_{r}\left(  \rho\right)  }=\chi\left(  \widehat
{S}_{r}\left(  \rho\right)  \right)  =2.\label{Gauss-Bonnett formul}%
\end{equation}

We need to estimate the Gauss-Bonnet integrand (see \cite{GW-82} p. 740ff or
\cite{ESS-89} p. 160ff). Define
\[
Q_{r}(\rho)\doteqdot G_{r}(\rho)-\det\left[  L_{r}(\rho)\right]  .
\]

>From Lemma \ref{lem smoothed set} (iii) we have the following estimate of the
sectional curvature $K_{x}$ of $\mathcal{M}$ at points in $\widehat{S}%
_{r}\left(  \rho\right)  $
\[
\sup_{x\in\widehat{S}_{r}\left(  \rho\right)  }K_{x}\leq\kappa(r+\rho-\eta
_{2})
\]
We estimate $Q_{r}(\rho)$ as on p. 160 of \cite{ESS-89}. Below $c(n)$ is a
constant depending only on $n$. For $r\geq1\geq\eta_{2}$
\begin{align*}
\left\vert Q_{r}(\rho)\right\vert  &  \leq c(n)\cdot\sum_{p=1}^{m}%
\kappa(r+\rho-\eta_{2})^{p}\cdot\left(  \left\vert L_{r}(\rho)\right\vert
^{2}\right)  ^{m-p}\\
&  \leq c(n)\cdot\sum_{p=1}^{m}\frac{a(r+\rho-\eta_{2})^{2p}}{(r+\rho-\eta
_{2})^{2p}}\cdot\left(  \frac{1}{\rho^{2}}\right)  ^{m-p}\text{ }\\
&  \leq c(n)\cdot\sum_{p=1}^{m}a(r+\rho-\eta_{2})^{2p}\cdot\frac{1}{\rho
^{n-1}}.
\end{align*}
It follows from the monotonicity of $a(r)$ that for $r\geq1$
\begin{equation}
\left\vert Q_{r}(\rho)\right\vert \leq c(n)\cdot\frac{1}{\rho^{n-1}}\cdot
\sum_{p=1}^{m}a(r-\eta_{2})^{2p}.\label{estimate for Q_r(rho)}%
\end{equation}

>From (\ref{estimate for L_r(rho)}), (\ref{Gauss-Bonnett formul}) and
(\ref{estimate for Q_r(rho)}) we get for $r\geq1$
\begin{align*}
\omega_{n-1} &  =\int_{\widehat{S}_{r}(\rho)}G_{r}(\rho)dV_{\widehat{S}%
_{r}\left(  \rho\right)  }\\
&  \leq\int_{\widehat{S}_{r}(\rho)}\left(  \left\vert Q_{r}(\rho)\right\vert
+\left\vert \det[L_{r}(\rho)]\right\vert \right)  dV_{\widehat{S}_{r}\left(
\rho\right)  }\\
&  \leq\int_{\widehat{S}_{r}(\rho)}\left(  c(n)\cdot\frac{1}{\rho^{n-1}}%
\cdot\sum_{p=1}^{m}a(r-\eta_{2})^{2p}+\frac{1}{\rho^{n-1}}\right)
dV_{\widehat{S}_{r}\left(  \rho\right)  }.
\end{align*}
Thus we have obtained the following lower bound for the areas of the hypersurfaces

\begin{lemma}
\label{estimate for area of S_r(rho)} \ There is a constant $c(n)$ depending
only on $n$ such that for any $r\geq1$ and $\rho\leq\rho(r-\eta_{2})=\frac
{\pi(r-\eta_{2})}{4a(r-\eta_{2})}$
\begin{equation}
\text{Area}\left[  \widehat{S}_{r}(\rho)\right]  \geq\frac{\omega_{n-1}%
}{c(n)\cdot\sum_{p=1}^{m}a(r-\eta_{2})^{2p}+1}\cdot\rho^{n-1}.
\label{vol S_r lower}%
\end{equation}

\end{lemma}

\subsection{The final argument\label{the argument}}

\begin{lemma}
\label{Lemma 20}(i) For $r\geq1$%
\[
\widehat{S}_{r}\subset\left.  B\left(  Q,\frac{r+\eta_{2}}{1-\theta\left(
r-\eta_{2}\right)  }\right)  \right\backslash B\left(  Q,r-\eta_{2}\right)  .
\]

(ii) For $r\geq1$ and $\rho\leq\rho(r-\eta_{2})$%
\[
\widehat{S}_{r}(\rho)\subset\left.  B\left(  Q,\frac{r+\eta_{2}}%
{1-\theta\left(  r-\eta_{2}\right)  }+\rho\right)  \right\backslash B\left(
Q,r-\eta_{2}\right)  .
\]

\end{lemma}

\begin{proof}
(i) For any $x\in\widehat{S}_{r}$ it follows from Lemma \ref{lem smoothed set}
(i) that $d(Q,x)\geq r-\eta_{2}$ and $b_{Q}(x)\leq r+\eta_{2}$. By Lemma
\ref{Busemann fcn sublevel containments}\ (ii), (i) we have for
\[
b_{Q}(x)\geq(1-\theta\left(  d(Q,x\right)  ))\cdot d(Q,x)\geq(1-\theta
(r-\eta_{2}))\cdot d(Q,x),
\]
Hence%
\[
d(Q,x)\leq\frac{r+\eta_{2}}{1-\theta\left(  r-\eta_{2}\right)  }.
\]

(ii) This follows from (i) and Lemma \ref{lem smoothed set} (iv).\medskip
\end{proof}

By Lemma \ref{Busemann fcn sublevel containments} (i) there is a $r_{2}%
=r_{2}(\mathcal{M})$ such that $\theta(r)<1/2$ for all $r\geq r_{2}$.

\begin{proposition}
\label{prop rel vol est below} Let $\left(  \mathcal{M}^{n},g\right)  $ be a
complete, noncompact Riemannian manifold with positive sectional curvature and
$Q\in\mathcal{M}$. For any $\eta_{1},\eta_{2}>0$ there is $r_{1}=r_{1}%
(\eta_{1},\mathcal{M})$ and $r_{2}=r_{2}(\mathcal{M})$ such that for any
$r\geq\max\left\{  r_{1},r_{2}\right\}  +5$
\begin{align*}
&  \text{Vol}\left[  B\left.  \left(  Q,\frac{r+\eta_{2}}{1-\theta\left(
r-\eta_{2}\right)  }+\rho(r-\eta_{2})\right)  \right\backslash B(Q,r-\eta
_{2})\right]  \\
&  \geq\frac{\omega_{n-1}}{c(n)\cdot\sum_{p=1}^{m}(\sqrt{ASCR(g)}+\eta
_{1})^{2p}+1}\cdot\frac{\rho(r-\eta_{2})^{n}}{n}.
\end{align*}

\end{proposition}

\begin{proof}
>From Lemma \ref{Lemma 20} (ii)%
\[
\cup_{0\leq\rho\leq\rho(r-\eta_{2})}\widehat{S}_{r}(\rho)\subset\left.
B\left(  Q,\frac{r+\eta_{2}}{1-\theta\left(  r-\eta_{2}\right)  }+\rho
(r-\eta_{2})\right)  \right\backslash B(Q,r-\eta_{2})
\]
we get by the area estimate (\ref{vol S_r lower})
\begin{align}
&  \text{Vol}\left[  B\left.  \left(  Q,\frac{r+\eta_{2}}{1-\theta\left(
r-\eta_{2}\right)  }+\rho(r-\eta_{2})\right)  \right\backslash B(Q,r-\eta
_{2})\right]  \nonumber\\
&  \geq\text{Vol}(\cup_{0\leq\rho\leq\rho(r-\eta_{2})}\widehat{S}_{r}%
(\rho))\nonumber\\
&  \geq\int_{0}^{\rho(r-\eta_{2})}\frac{\omega_{n-1}}{c(n)\cdot\sum_{p=1}%
^{m}a(r-\eta_{2})^{2p}+1}\cdot\rho^{n-1}d\rho\nonumber\\
&  =\frac{\omega_{n-1}}{c(n)\cdot\sum_{p=1}^{m}a(r-\eta_{2})^{2p}+1}\cdot
\frac{\rho(r-\eta_{2})^{n}}{n}.\label{vol for ball from hypsurf}%
\end{align}
The proposition follows from $\theta\left(  r-\eta_{2}\right)  \leq\frac{1}%
{2}$ and $a(r-\eta_{2})\leq\sqrt{ASCR(g)}+\eta_{1}$ when $r\geq r\geq
\max\left\{  r_{1},r_{2}\right\}  +5$.\medskip
\end{proof}

Now we finish the proof of Theorem \ref{Main-ASCR}. Let $\left(
\varepsilon_{0},k_{0},L_{0}\right)  =\left(  \min(\varepsilon_{a}%
,\varepsilon_{b}),k_{a},L_{b}\right)  $ given in Proposition
\ref{prop rel vol comp}. Let $N:S^{n-1}\times\lbrack-L,L]\rightarrow
\mathcal{M}$ be a $\left(  \varepsilon,k,L\right)  $-neck with $\left(
\varepsilon,k,L\right)  \preccurlyeq\left(  \varepsilon_{0},k_{0}%
,L_{0}\right)  $ and $Q\in N(S^{n-1}\times\{-L_{0}\})$; this will give us a
constant $r_{0}$. Choose positive $\eta_{1}<1$\ and $\eta_{2}<1$, this will
give two constants $r_{1}$ and $r_{2}$. Choosing an $r\geq\max\left\{
r_{0},r_{1},r_{2}\right\}  +5$ and letting $R_{2}\doteqdot\frac{r+\eta_{2}%
}{1-\theta\left(  r-\eta_{2}\right)  }+\rho(r-\eta_{2})$ and $R_{1}\doteqdot
r-\eta_{2}$ in Proposition \ref{prop rel vol comp}, we get from
(\ref{eq rel vol up est}) and (\ref{vol for ball from hypsurf})%

\begin{align*}
&  \frac{\left(  \frac{11}{10}\right)  ^{n}\cdot\omega_{n-1}\cdot\frac{10}%
{9}\cdot12}{\frac{\omega_{n-1}}{n}\cdot((\frac{9}{10}L_{0}+2)^{n}-(\frac
{9}{10}L_{0}-2)^{n})}\\
&  \geq\frac{\frac{\omega_{n-1}}{c(n)\cdot\sum_{p=1}^{m}a(r-\eta_{2})^{2p}%
+1}\cdot\frac{\rho(r-\eta_{2})^{n}}{n}}{\frac{\omega_{n-1}}{n}\cdot\left(
\left[  \frac{r+\eta_{2}}{1-\theta\left(  r-\eta_{2}\right)  }+\rho(r-\eta
_{2})\right]  ^{n}-(r-\eta_{2})^{n}\right)  }\\
&  \geq\frac{1}{c(n)\cdot\sum_{p=1}^{m}a(r-\eta_{2})^{2p}+1}\cdot\frac
{1}{\left(  \frac{r+\eta_{2}}{\left[  1-\theta\left(  r-\eta_{2}\right)
\right]  \cdot\rho(r-\eta_{2})}+1\right)  ^{n}}\\
&  \geq\frac{1}{c(n)\cdot\sum_{p=1}^{m}a(r-\eta_{2})^{2p}+1}\cdot\frac
{1}{\left(  \frac{12}{\pi}a(r-\eta_{2})+1\right)  ^{n}},
\end{align*}
where in the last inequality we have used $\frac{r+\eta_{2}}{\rho(r-\eta_{2}%
)}\leq\frac{6}{\pi}a(r-\eta_{2})$ and $\theta\left(  r-\eta_{2}\right)
\leq1/2.$ This is because $r\geq5$ and $\frac{r-\eta_{2}}{\rho(r-\eta_{2}%
)}=\frac{4}{\pi}a(r-\eta_{2})$.\label{changed leq to =} Using $a(r-\eta
_{2})\leq\sqrt{ASCR(g)}+\eta_{1}$ we get
\begin{align*}
&  \frac{\left(  \frac{11}{10}\right)  ^{n}\cdot\omega_{n-1}\cdot1^{n}%
\cdot\frac{10}{9}\cdot12}{\frac{\omega_{n-1}}{n}\left(  (\frac{9}{10}%
L_{0}+2)^{n}-(\frac{9}{10}L_{0}-2)^{n}\right)  }\\
&  \geq\frac{1}{c(n)\cdot\sum_{p=1}^{m}\left(  \sqrt{ASCR(g)}+\eta_{1}\right)
^{2p}+1}\cdot\frac{1}{\left(  \frac{12}{\pi}\left(  \sqrt{ASCR(g)}+\eta
_{1}\right)  +1\right)  ^{n}}.
\end{align*}
Now it is clear that we can make $ASCR(g)\geq C_{0}$ if we choose $L_{0}$
large enough. Theorem \ref{Main-ASCR} is proved.\textrm{ }

\section{Existence of necklike points in ancient
solutions\label{Existence of necklike points in ancient solutions}}

We first recall the following result of Hamilton which is based on an estimate
of Hamilton (Theorem 24.4 in \cite{H-95a}) and Ivey \cite{I-93}.

\begin{lemma}
\label{Ancient3ManifoldsHaveNonnegativeSectional}If $\left(  \mathcal{M}%
^{3},g\left(  t\right)  \right)  ,\ t\in\left(  -\infty,\omega\right)  ,$ is a
complete ancient solution of the Ricci flow with bounded curvature, then
$g\left(  t\right)  $ has nonnegative sectional curvature for all $t\in\left(
-\infty,\omega\right)  .$
\end{lemma}

\begin{proof}
See \cite{H-95a} or for a more detailed proof \cite{CKE}.
\end{proof}

We also note that in a lemma in \S 19 of \cite{H-95a}, Hamilton also proved
that if $\left(  \mathcal{M}^{n},g\left(  t\right)  \right)  $ is a non-Ricci
flat ancient solution with nonnegative Ricci curvature, then there exists a
constant $c_{0}>0$ depending only on $n$ such that
\[
\liminf_{t\rightarrow-\infty}\left\vert t\right\vert \sup_{x\in\mathcal{M}%
}R\left(  x,t\right)  \geq c_{0}.
\]
(See also Lemma 19.4 of \cite{H-95a} for a related result for Type I
singularities.) In some sense this may be considered as an elementary gap-type
result for ancient solutions to the Ricci flow regarding the geometry at
$t=-\infty$ of the space-time manifold $\mathcal{M}^{n}\times\left(
-\infty,\omega\right)  .$

Following the terminology of \cite{CKE} we shall say that $\left(  x,t\right)
$ is an \textit{ancient Type I-like }$c$\textit{-essential point} if
\[
\left\vert \operatorname*{Rm}\left(  x,t\right)  \right\vert \cdot\left\vert
t\right\vert \geq c>0.
\]
We say that $\left(  x,t\right)  $ is a $\delta$\textit{-necklike point} if
there exists a unit $2$-form $\theta$ at $\left(  x,t\right)  $ such that
\[
\left\vert \operatorname*{Rm}-R\left(  \theta\otimes\theta\right)  \right\vert
\leq\delta\left\vert \operatorname*{Rm}\right\vert .
\]
We have the following result due to the first named author which we believe to
have first appeared in the unpublished notes \cite{CKE}. Similar results were
proved earlier by Hamilton in Theorem 24.6 of \cite{H-95a} on the existence of
necklike points in Type I singular solutions to the Ricci flow on compact
$3$-manifolds on finite time intervals $[0,T)$ (see also Theorem 3.3 in \S 2.3
of \cite{H-97} for a similar result in dimension four) and Corollary 3.5 of
\S 2.3 of \cite{H-97} on the existence of necklike points in Type I-like
ancient solutions with bounded positive isotropic curvature.

\begin{theorem}
\label{Type I like necklike points in ancient solutions}Let $\left(
\mathcal{M}^{3},g\left(  t\right)  \right)  ,$ $t\in\left(  -\infty
,\omega\right)  ,$ be a complete ancient solution of the Ricci flow with
bounded positive sectional curvature. Suppose that
\[
\sup_{\mathcal{M}\times(-\infty,0]}\left|  t\right|  ^{\gamma}R\left(
x,t\right)  <\infty
\]
for some $\gamma>0$. Then either

\begin{enumerate}
\item $\left(  \mathcal{M},g\left(  t\right)  \right)  $ is isometric to a
shrinking spherical space form, or

\item there exists a constant $c>0$ such that for all $\tau\in(-\infty,0]$ and
$\delta>0$, there exist $x\in\mathcal{M}$ and $t\in\left(  -\infty
,\tau\right)  $ such that $\left(  x,t\right)  $ is an ancient Type I-like
$c$-essential point and a $\delta$-necklike point.
\end{enumerate}
\end{theorem}

\begin{proof}
We shall show that if for every $c>0$ there exist $\tau\in(-\infty,0]$ and
$\delta>0$ such that there are no ancient Type I-like $c$-essential $\delta
$-necklike points before time $\tau$, then $(\mathcal{M},g(t))$ is isometric
to a shrinking spherical space form.

By the hypothesis, there exists $\gamma>0$ such that
\[
K\doteqdot\sup_{\mathcal{M}\times(-\infty,0]}\left|  t\right|  ^{\gamma
}R\left(  x,t\right)  <\infty.
\]
(When $\gamma=1$, this is the definition of an ancient Type I-like solution.)
Since the scalar curvature of $(\mathcal{M},g(t))$ is positive, the function
\[
G\doteqdot\left|  t\right|  ^{\gamma\varepsilon/2}\frac{|\overset{\circ
}{\operatorname*{Rm}}|^{2}}{R^{2-\varepsilon}}%
\]
is well-defined. Since the sectional curvatures are positive, $|\overset
{\circ}{\operatorname*{Rm}}|\leq\left|  \operatorname*{Rm}\right|  \leq R,$
and we have the estimates
\[
G\leq R^{\varepsilon}\left|  t\right|  ^{\gamma\varepsilon/2}\leq
K^{\varepsilon}\left|  t\right|  ^{-\gamma\varepsilon/2},
\]
which show that $G$ is bounded for all times $-\infty<t\leq0$ and satisfies
\begin{equation}
\lim_{t\rightarrow-\infty}\sup_{x\in\mathcal{M}}G\left(  x,t\right)  =0.
\label{GDecaysAtAncientTimes}%
\end{equation}
If there are no ancient $c$-essential $\delta$-necklike points on the time
interval $(-\infty,\tau]$, then for every $x\in\mathcal{M}$ and $t\in\left(
-\infty,\tau\right)  $ either
\begin{equation}
\left|  \operatorname*{Rm}\left(  x,t\right)  \right|  \cdot\left|  t\right|
<c \label{NoAncientcEssential}%
\end{equation}
or we have
\begin{equation}
\left|  \operatorname*{Rm}-R\left(  \theta\otimes\theta\right)  \right|
>\delta\left|  \operatorname*{Rm}\right|  \label{NoAncientDeltaNecklike}%
\end{equation}
for every unit $2$-form $\theta$ at $\left(  x,t\right)  $.

A straightforward computation yields that if $\phi$ is a nonnegative function
and $\psi$ is a positive function, both defined on space and time, then
\begin{align*}
\left(  \frac{\partial}{\partial t}-\Delta\right)  \left(  \frac{\phi^{\alpha
}}{\psi^{\beta}}\right)   &  =\alpha\frac{\phi^{\alpha-1}}{\psi^{\beta}%
}\left(  \frac{\partial}{\partial t}-\Delta\right)  \phi-\beta\frac
{\phi^{\alpha}}{\psi^{\beta+1}}\left(  \frac{\partial}{\partial t}%
-\Delta\right)  \psi\\
&  -\alpha\left(  \alpha-1\right)  \frac{\phi^{\alpha-2}}{\psi^{\beta}%
}\left\vert \nabla\phi\right\vert ^{2}-\beta\left(  \beta+1\right)  \frac
{\phi^{\alpha}}{\psi^{\beta+2}}\left\vert \nabla\psi\right\vert ^{2}\\
&  +2\alpha\beta\frac{\phi^{\alpha-1}}{\psi^{\beta+1}}\left\langle \nabla
\phi,\nabla\psi\right\rangle .
\end{align*}
Taking $\phi=\left(  -t\right)  ^{\gamma\varepsilon/2}|\overset{\circ
}{\operatorname*{Rm}}|^{2}$, $\psi=R$, $\alpha=1$, and $\beta=2-\varepsilon,$
a computation yields
\[
\frac{\partial}{\partial t}G\leq\Delta G+\frac{2\left(  1-\varepsilon\right)
}{R}\left\langle \nabla G,\nabla R\right\rangle +2J,
\]
where
\[
J\doteqdot\frac{\left\vert t\right\vert ^{\gamma\varepsilon/2}}%
{R^{3-\varepsilon}}\left[  \varepsilon|\overset{\circ}{\operatorname*{Rm}%
}|^{2}\left(  \left\vert \operatorname*{Rm}\right\vert ^{2}-\frac{\gamma
R}{4\left\vert t\right\vert }\right)  -P\right]
\]
and
\[
P\doteqdot\lambda^{2}\left(  \mu-\nu\right)  ^{2}+\mu^{2}\left(  \lambda
-\nu\right)  ^{2}+\nu^{2}\left(  \lambda-\mu\right)  ^{2}\geq0
\]
(here $\lambda,\mu,\nu$ are the eigenvalues of $Rm.$) Fix any $\left(
x,t\right)  $ with $t<\tau\leq0$. If the first alternative
(\ref{NoAncientcEssential}) holds there with $c\leq\gamma/8$, then we may
estimate the term
\[
\left\vert \operatorname*{Rm}\right\vert ^{2}-\frac{\gamma R}{4\left\vert
t\right\vert }\leq R\left(  \left\vert \operatorname*{Rm}\right\vert
-\frac{\gamma}{4\left\vert t\right\vert }\right)  <R\left(  \frac
{c}{\left\vert t\right\vert }-\frac{\gamma}{4\left\vert t\right\vert }\right)
\leq-\frac{\gamma R}{8\left\vert t\right\vert }%
\]
and hence dropping the $-P\leq0$ term yields
\[
J\leq-\frac{\gamma\varepsilon}{8\left\vert t\right\vert }G.
\]
On the other hand if the second alternative (\ref{NoAncientDeltaNecklike})
holds, then Lemma \ref{NotProductlike} in Appendix C implies there exists
$\eta\left(  \delta\right)  >0$ such that $P\geq\eta\left\vert
\operatorname*{Rm}\right\vert ^{2}|\overset{\circ}{\operatorname*{Rm}}|^{2}$
whence it follows that taking $\varepsilon\leq\eta$ gives
\[
J\leq-\frac{\gamma\varepsilon}{4\left\vert t\right\vert }G.
\]
Thus in either case, $G$ is a subsolution of the heat equation for all times
$-\infty<t<\tau$, because
\[
\frac{\partial}{\partial t}G\leq\Delta G+\frac{2\left(  1-\varepsilon\right)
}{R}\left\langle \nabla G,\nabla R\right\rangle -\frac{\gamma\varepsilon
}{2\left\vert t\right\vert }G.
\]
By the weak maximum principle, which applies even when $\mathcal{M}$ is
noncompact since the both the curvatures and $G$ are bounded, $\sup
_{x\in\mathcal{M}}G\left(  x,t\right)  $ is a nonincreasing function of time.
We may then use (\ref{GDecaysAtAncientTimes}) to conclude that $G\equiv0$,
hence that $\left(  \mathcal{M}^{3},g(t)\right)  $ locally isometric to a
round $\mathcal{S}^{3}$. Since $\left(  \mathcal{M}^{3},g(t)\right)  $ is
complete, we conclude that it is compact and globally isometric to a spherical
space form $\mathcal{S}^{3}/\Gamma$.
\end{proof}

Now we shall assume that $\mathcal{M}^{3}$ is noncompact. Since $g\left(
t\right)  $ has positive sectional curvature, by a result of Gromoll and Meyer
\cite{GM-69}, there is an injectivity radius estimate:

\begin{proposition}
$\mathcal{M}^{3}$ is diffeomorphic to $\mathbb{R}^{3}$ and%
\[
\text{inj}\left(  \mathcal{M}^{3},g\left(  t\right)  \right)  \geq\frac{\pi
}{\sqrt{K_{\sup}\left(  t\right)  }}%
\]
where $K_{\sup}\left(  t\right)  $ is the supremum of the sectional curvatures
of $g\left(  t\right)  .$
\end{proposition}

We may apply this estimate and a standard compactness theorem to obtain a
cylinder limit solution. In particular, we have:

\begin{theorem}
\label{existence of cylinder limit}If $\left(  \mathcal{M}^{3},g\left(
t\right)  \right)  ,$ $t\in\left(  -\infty,\omega\right)  ,$ is a complete
noncompact Type I-like ancient solution of the Ricci flow with bounded
positive sectional curvature on an orientable $3$-manifold, then there exists
a sequence of points and times $\left(  x_{i},t_{i}\right)  \in\mathcal{M}%
^{3}\times\left(  -\infty,\omega\right)  $ such that the dilated and
translated solutions $\left(  \mathcal{M}^{3},g_{i}\left(  t\right)
,x_{i}\right)  ,$ $t\in\left(  -\infty,\omega_{i}\right)  ,$ where%
\[
g_{i}\left(  t\right)  =R\left(  x_{i},t_{i}\right)  \cdot g\left(
t_{i}+\frac{t}{R\left(  x_{i},t_{i}\right)  }\right)
\]
and $\omega_{i}=R\left(  x_{i},t_{i}\right)  \left(  T-t_{i}\right)  ,$ limit
to a solution $\left(  \mathcal{M}_{\infty}^{3},g_{\infty}\left(  t\right)
,x_{\infty}\right)  ,$ $t\in\left(  -\infty,\omega_{\infty}\right)  ,$ to the
Ricci flow isometric to the standard shrinking cylinder $S^{2}\times
\mathbb{R}.$
\end{theorem}

\begin{remark}
The conclusion also holds with the additional condition on the sequence of
times: $\lim_{i\rightarrow\infty}t_{i}=-\infty.$
\end{remark}

\begin{proof}
Choose any sequence $\left\{  \delta_{i}\right\}  _{i\in\mathbb{N}}$ with
$\lim_{i\rightarrow\infty}\delta_{i}=0.$ By Theorem
\ref{Type I like necklike points in ancient solutions} (with $\gamma=1,$)
since $\mathcal{M}$ is noncompact\ (and hence $\left(  \mathcal{M},g\left(
t\right)  \right)  $ is not isometric to a shrinking spherical space form),
there exists a constant $c>0$ such that there exists points and times $\left(
x_{i},t_{i}\right)  $ such that
\[
\left\vert \operatorname*{Rm}\left(  x_{i},t_{i}\right)  \right\vert
\cdot\left\vert t_{i}\right\vert \geq c>0.
\]
and%
\[
\left\vert \operatorname*{Rm}-R\left(  \theta\otimes\theta\right)  \right\vert
\left(  x_{i},t_{i}\right)  \leq\delta_{i}\left\vert \operatorname*{Rm}\left(
x_{i},t_{i}\right)  \right\vert .
\]
for some unit $2$-forms $\theta_{i}.$ By Gromoll and Meyer's injectivity
radius estimate, we may apply Hamilton's Gromov-type compactness theorem for
solutions of the Ricci flow \cite{H-95b} to the sequence of pointed solutions
$\left(  \mathcal{M}^{3},g_{i}\left(  t\right)  ,x_{i}\right)  ,$ $t\in\left(
-\infty,\omega_{i}\right)  $ where%
\[
g_{i}\left(  t\right)  =R\left(  x_{i},t_{i}\right)  \cdot g\left(
t_{i}+\frac{t}{R\left(  x_{i},t_{i}\right)  }\right)
\]
and $\omega_{i}=R\left(  x_{i},t_{i}\right)  \left(  T-t_{i}\right)  .$ We
obtain a complete limit ancient solution to the Ricci flow $\left(
\mathcal{M}_{\infty}^{3},g_{\infty}\left(  t\right)  ,x_{\infty}\right)  ,$
$t\in\left(  -\infty,\omega_{\infty}\right)  ,$ on a noncompact\footnote{If
$\mathcal{M}_{\infty}^{3}$ were compact, then $\mathcal{M}^{3}$ would be
diffeomorphic to $\mathcal{M}_{\infty}^{3},$ which contradicts the asumption
that $\mathcal{M}^{3}$ is noncompact.} orientable $3$-manifold with bounded
nonnegative sectional curvature and%
\[
\operatorname*{Rm}\left(  g_{\infty}\left(  x_{\infty},0\right)  \right)
=R\left(  g_{\infty}\left(  x_{\infty},0\right)  \right)  \left(
\theta_{\infty}\otimes\theta_{\infty}\right)
\]
for some unit $2$-form $\theta_{\infty}.$ By the strong maximum principle we
conclude that the universal covering solution $\left(  \widetilde{\mathcal{M}%
}_{\infty}^{3},\tilde{g}_{\infty}\left(  t\right)  \right)  ,$ $t\in\left(
-\infty,\omega_{\infty}\right)  ,$ is isometric to the product of $\mathbb{R}$
and a complete ancient Type I-like solution $\left(  \mathcal{N}^{2}%
,h_{\infty}\left(  t\right)  \right)  ,$ $t\in\left(  -\infty,\omega_{\infty
}\right)  ,$ to the Ricci flow on a surface with bounded positive curvature.
By Theorem 26.1 of \cite{H-95a}, $\left(  \mathcal{N}^{2},h_{\infty}\left(
t\right)  \right)  $ is isometric to a shrinking round sphere (and in
particular, $\mathcal{N}^{2}$ is compact.) Now there are only two noncompact
orientable quotients of $S^{2}\times\mathbb{R}$ : $S^{2}\times\mathbb{R}$
itself and $\mathbb{R}\tilde{\times}\mathbb{R}P^{2},$ the nontrivial
$\mathbb{R}$-bundle over $\mathbb{R}P^{2}.$ Topologically, $\mathbb{R}%
\tilde{\times}\mathbb{R}P^{2}\cong S^{2}\times\lbrack0,\infty)/\sim$ where
$\left(  x,0\right)  \sim\left(  -x,0\right)  ,$ which is diffeomorphic to
$\mathbb{R}P^{3}-\bar{B}^{3}$. If $\mathcal{M}_{\infty}^{3}\cong
\mathbb{R}\tilde{\times}\mathbb{R}P^{2},$ then $\mathcal{M}^{3}$ admits an
embedded (one-sided) $\mathbb{R}P^{2}$ (since $\mathcal{M}_{\infty}^{3}$
does). However, by the work of Gromoll and Meyer \cite{GM-69}, we know that
$\mathcal{M}^{3}$ is diffeomorphic to $\mathbb{R}^{3}.$ This yields a
contradiction. Hence $\left(  \mathcal{M}_{\infty}^{3},g_{\infty}\left(
t\right)  ,x_{\infty}\right)  $ is isometric to $\mathbb{R}\times\left(
\mathcal{N}^{2},h_{\infty}\left(  t\right)  \right)  ,$ where $\left(
\mathcal{N}^{2},h_{\infty}\left(  t\right)  \right)  $ is a shrinking round
$2$-sphere.\medskip
\end{proof}

We are now in a position to prove Theorem
\ref{Type I ancient have infinite ASCR} assuming Theorem \ref{Main-ASCR}%
.\medskip

\begin{proof}
[Proof of Theorem \ref{Type I ancient have infinite ASCR}]Let $\left(
\mathcal{M}^{3},g\left(  t\right)  \right)  ,$ $-\infty<t<\omega,$ be a
complete orientable noncompact ancient Type I-like solution to the Ricci flow
with bounded positive sectional curvature. By Theorem
\ref{existence of cylinder limit}, there exists a sequence of points and times
$\left(  x_{i},t_{i}\right)  \in\mathcal{M}^{3}\times\left(  -\infty
,\omega\right)  $ such that the dilated and translated solutions $\left(
\mathcal{M}^{3},g_{i}\left(  t\right)  ,x_{i}\right)  ,$ $t\in\left(
-\infty,\omega_{i}\right)  ,$ limit to the standard shrinking cylinder
$S^{2}\times\mathbb{R}.$ Now, by Theorem \ref{Main-ASCR}, given any
$A_{0}<\infty,$ there exists $\left(  \varepsilon_{0},k_{0},L_{0}\right)  $
such that if $\left(  \mathcal{M}^{3},g\right)  $ is a complete, noncompact
Riemannian manifold with bounded positive sectional curvature and containing
an $\left(  \varepsilon,k,L\right)  $-neck with $\left(  \varepsilon
,k,L\right)  \preccurlyeq\left(  \varepsilon_{0},k_{0},L_{0}\right)  ,$ then
~ASCR$\left(  g\right)  \geq A_{0}.$ On the other hand, since $\left(
\mathcal{M}^{3},g_{i}\left(  0\right)  ,x_{i}\right)  $ limits to a standard
cylinder $S^{2}\times\mathbb{R},$ for $i$ large enough, there exists an
$\left(  \varepsilon_{0},k_{0},L_{0}\right)  $-neck in $\left(  \mathcal{M}%
^{3},g_{i}\left(  0\right)  ,x_{i}\right)  .$ This implies
\[
\text{ASCR}\left(  g\left(  t_{i}\right)  \right)  =\text{ASCR}\left(
g_{i}\left(  0\right)  \right)  \geq A_{0}.
\]
In \S 19 of \cite{H-95a} Hamilton has proved that the asymptotic scalar
curvature ratio of a complete ancient solution to the Ricci flow with bounded
nonnegative curvature operator (when $n=3,$ this is the same as nonnegative
sectional curvature) is constant in time. This implies that for all
$t\in\left(  -\infty,\omega\right)  $ we have ASCR$\left(  g\left(  t\right)
\right)  \geq A_{0}.$ Since $A_{0}<\infty$ is arbitrary, we conclude that
\[
ASCR\left(  g\left(  t\right)  \right)  \equiv\infty
\]
for all $t\in\left(  -\infty,\omega\right)  .$
\end{proof}

\section{Appendices}

\subsection{A}

\begin{proof}
[Proof of Lemma \ref{lem absolute neck}]It suffices to show that there is
$\varepsilon^{\prime}$ such that if $N:S^{n-1}\times\left[  a,b\right]
\rightarrow\mathcal{M}$ is an $\left(  \varepsilon^{\prime},k,L\right)  $-neck
then inequalities (\ref{absol neck ineq 1}) and (\ref{absol neck ineq 2}) are
satisfied on $z\in\lbrack\frac{a+b}{2}-L,\frac{a+b}{2}+L]$.

Assume $\varepsilon^{\prime}\leq\min\left\{  1,\frac{\ln2}{2L}\right\}  $ so
that $\left\vert \frac{d\ln r(z)}{dz}\right\vert \leq\varepsilon^{\prime}.$
Using the inequality $|e^{x}-1|\leq2|x|$ for $|x|\leq\ln2,$ we get for
$z\in\lbrack\frac{a+b}{2}-L,\frac{a+b}{2}+L],$
\[
\left\vert \ln\frac{r^{2}(z)}{r^{2}(\frac{a+b}{2})}\right\vert \leq2\left\vert
\int_{\left(  a+b\right)  /2}^{z}\frac{d}{d\zeta}\left(  \ln r\left(
\zeta\right)  \right)  d\zeta\right\vert \leq2L\varepsilon^{\prime}\leq\ln2,
\]
which implies%
\[
\left\vert \frac{r^{2}(z)}{r^{2}(\frac{a+b}{2})}-1\right\vert \leq2\left\vert
\ln\frac{r^{2}(z)}{r^{2}(\frac{a+b}{2})}\right\vert \leq4\varepsilon^{\prime
}L.
\]
>From $\left\vert \frac{1}{r^{2}(z)}g-\bar{g}\right\vert _{\bar{g}}%
\leq\varepsilon^{\prime}$ we get%
\[
\left\vert \frac{1}{r^{2}(z)}\cdot g\right\vert _{\bar{g}}\leq\sqrt
{n}+\varepsilon^{\prime}%
\]
since $\left\vert \bar{g}\right\vert _{\bar{g}}^{2}=\bar{g}^{ij}\bar{g}%
_{ij}=n.$ Hence for $z\in\lbrack\frac{a+b}{2}-L,\frac{a+b}{2}+L]$%
\begin{align*}
\left\vert \frac{1}{r^{2}(\frac{a+b}{2})}\cdot g-\bar{g}\right\vert _{\bar
{g}}  &  \leq\left\vert \frac{1}{r^{2}(z)}\cdot g-\bar{g}\right\vert _{\bar
{g}}+\left\vert \left(  \frac{r^{2}(z)}{r^{2}(\frac{a+b}{2})}-1\right)
\cdot\left(  \frac{1}{r^{2}(z)}\cdot g\right)  \right\vert _{\bar{g}}\\
&  \leq\varepsilon^{\prime}+4\varepsilon^{\prime}L\cdot(\sqrt{n}%
+\varepsilon^{\prime}).
\end{align*}
If we choose
\[
\varepsilon^{\prime}\leq\min\left\{  1,\frac{\ln2}{2L},\frac{\varepsilon
}{1+8L\sqrt{n}}\right\}  ,
\]
inequality (\ref{absol neck ineq 1}) will hold for $z\in\lbrack\frac{a+b}%
{2}-L,\frac{a+b}{2}+L]$.

Let
\[
R(z)=\ln\frac{r^{2}(z)}{r^{2}(\frac{a+b}{2})}.
\]
Assume now $\varepsilon^{\prime}\leq\min\left\{  \frac{1}{2},\frac{\ln2}%
{2L}\right\}  $. From $\left\vert \frac{d^{j}\log r\left(  z\right)  }{dz^{j}%
}\right\vert \leq\varepsilon^{\prime}$ for $1\leq j\leq k$ and $e^{R(z)}%
\leq1+4\varepsilon^{\prime}L$ we get for $1\leq j\leq k$ and $z\in\lbrack
\frac{a+b}{2}-L,\frac{a+b}{2}+L]$%
\[
\left\vert \frac{d^{j}R(z)}{dz^{j}}\right\vert \leq2\varepsilon^{\prime}%
\]
and%
\begin{align*}
\left\vert \frac{d^{j}}{dz^{j}}\frac{r^{2}(z)}{r^{2}(\frac{a+b}{2}%
)}\right\vert  &  =\left\vert \frac{d^{j}}{dz^{j}}e^{R(z)}\right\vert \leq
c(j)\cdot e^{R(z)}\cdot\max_{i_{1}+\cdots+i_{p}=j}\left\{  \left\vert
\frac{d^{i_{1}}R(z)}{dz^{i_{1}}}\bullet\cdots\bullet\frac{d^{i_{p}}%
R(z)}{dz^{i_{p}}}\right\vert \right\}  \\
&  \leq c(j)\cdot(1+4\varepsilon^{\prime}L)\cdot2\varepsilon^{\prime}%
\end{align*}
where $c(j)$ is a constant depending only on $j$. Hence for $1\leq j\leq k$
and $z\in\lbrack\frac{a+b}{2}-L,\frac{a+b}{2}+L]$%
\begin{align*}
&  \left\vert \bar{\nabla}^{j}\left(  \frac{1}{r^{2}\left(  \frac{a+b}%
{2}\right)  }\cdot g\right)  \right\vert _{\bar{g}}=\left\vert \bar{\nabla
}^{j}\left(  e^{R(z)}\cdot\frac{1}{r^{2}\left(  z\right)  }\cdot g\right)
\right\vert _{\bar{g}}\\
&  \leq\sum_{i=0}^{j}\left(
\begin{array}
[c]{c}%
j\\
i
\end{array}
\right)  \cdot\left\vert \frac{d^{i}}{dz^{i}}e^{R(z)}\right\vert
\cdot\left\vert \bar{\nabla}^{j-i}\left(  \frac{1}{r^{2}\left(  z\right)
}\cdot g\right)  \right\vert _{\bar{g}}\\
&  \leq(1+4\varepsilon^{\prime}L)\varepsilon^{\prime}+\sum_{i=1}^{j-1}\left(
\begin{array}
[c]{c}%
j\\
i
\end{array}
\right)  \cdot c(i)\cdot(1+4\varepsilon^{\prime}L)\cdot2\varepsilon^{\prime
}\cdot\varepsilon^{\prime}+c(j)\cdot(1+4\varepsilon^{\prime}L)\cdot
2\varepsilon^{\prime}.
\end{align*}
Note that in the last inequality above we have used $\left\vert \bar{\nabla
}^{j-i}\left(  \frac{1}{r^{2}\left(  z\right)  }\cdot g\right)  \right\vert
_{\bar{g}}\leq\varepsilon^{\prime}$ for $1\leq j-i\leq k$. So there is a
constant $C_{k}$ depending only on $k$, such that if we choose
\[
\varepsilon^{\prime}\leq\min\left\{  \frac{1}{2},\frac{\ln2}{2L}%
,\frac{\varepsilon}{C_{k}\cdot(1+4L)}\right\}  ,
\]
then inequality (\ref{absol neck ineq 2}) will hold for $z\in\lbrack\frac
{a+b}{2}-L,\frac{a+b}{2}+L]$. The lemma is proved.\medskip
\end{proof}

\subsection{B}

\begin{proof}
[Proof of Lemma \ref{lem smoothed set}]We shall use the properties of Busemann
functions in Lemma \ref{lem prop busemann}.

\textbf{(i)} For any $x\in\widehat{S}_{r},\widehat{b}_{Q}(x)=r$. It follows
from (\ref{eq for approx busemann}) that $|b_{Q}(x)-r|<\eta_{2}$. Hence
$\widehat{S}_{r}\subset b_{Q}^{-1}((r-\eta_{2},r+\eta_{2}))$.

For any $x\in\widehat{S}_{r}$, $b_{Q}(x)>r-\eta_{2}.$ Since $b_{Q}(x)\leq
d(Q,x)$, $d(Q,x)>$ $r-\eta_{2}$ and $\widehat{S}_{r}\subset\mathcal{M}%
\backslash\overline{B}(Q,r-\eta_{2})$.

\textbf{(ii)} For any $x\in B(Q,r-\eta_{2})$, $d(Q,x)<r-\eta_{2}$. Since
$b_{Q}(x)\leq d(Q,x)$, $b_{Q}(x)<r-\eta_{2}$ and $B(Q,r-\eta_{2})\subset
b_{Q}^{-1}((-\infty,r-\eta_{2}))$.

For any $x\in b_{Q}^{-1}((-\infty,r-\eta_{2}))$, $b_{Q}(x)<r-\eta_{2}$. By
(\ref{eq for approx busemann}), $\widehat{b}_{Q}(x)<b_{Q}(x)+\eta_{2}<r$.
Hence $b_{Q}^{-1}((-\infty,r-\eta_{2}))\subset\widehat{C}_{r}$.

For any $x\in\widehat{C}_{r}$, $\widehat{b}_{Q}(x)<r$. By
(\ref{eq for approx busemann}), $b_{Q}(x)<\widehat{b}_{Q}(x)+\eta_{2}%
<r+\eta_{2}$. Hence $\widehat{C}_{r}\subset b_{Q}^{-1}((-\infty,r+\eta_{2})).$

\textbf{(iii)} For any $x\in\widehat{S}_{r}(\rho)$, $d(x,\widehat{S}%
_{r})=d(x,\widehat{C}_{r})=\rho$. Let $\gamma$ be a minimal geodesic from $Q$
to $x$. It is clear that $\gamma$ intersects $\widehat{S}_{r}$ at some point
$\gamma(l)$, thus%
\[
d(Q,x)=d(Q,\gamma(l))+d(\gamma(l),x)\geq d(Q,\widehat{S}_{r})+d(\widehat
{S}_{r},x)>r-\eta_{2}+\rho,
\]
and hence $\widehat{S}_{r}(\rho)\subset\mathcal{M}\backslash\overline
{B}(Q,r+\rho-\eta_{2})$.

\textbf{(iv)} For any $x\in\widehat{S}_{r}(\rho)$ let $y\in\widehat{S}_{r}$ be
a point satisfying $d(x,y)=d(x,\widehat{S}_{r})=\rho$. Then from $\widehat
{S}_{r}\subset B(Q,\eta)$
\[
d(Q,x)\leq d(Q,y)+d(y,x)<\eta+\rho,
\]
so that $\widehat{S}_{r}(\rho)\subset B(Q,\eta+\rho)$.
\end{proof}

\subsection{C}

Here we give the proof of an estimate of Hamilton used in the proof of Theorem
\ref{Type I like necklike points in ancient solutions}.\footnote{This is a
version of an estimate in the proof of Theorem 24.6 of \cite{H-95a}.}

\begin{lemma}
\label{NotProductlike}If for some $\delta\in\left(  0,1\right)  $ we have
$\left|  \operatorname*{Rm}-R\left(  \theta\otimes\theta\right)  \right|
^{2}\geq\delta\left|  \operatorname*{Rm}\right|  ^{2}$ for every unit $2$-form
$\theta$, then
\[
P\geq\frac{\delta}{96\left(  3-\delta\right)  }\left|  \operatorname*{Rm}%
\right|  ^{2}|\overset{\circ}{\operatorname*{Rm}}|^{2}.
\]

\end{lemma}

\begin{proof}
We may assume without loss of generality that $\left\vert \lambda\right\vert
\geq\left\vert \mu\right\vert \geq\left\vert \nu\right\vert $. The hypothesis
implies that
\[
\mu^{2}+\nu^{2}+\mu\nu\geq\frac{\delta}{2}\left(  \lambda^{2}+\mu^{2}+\nu
^{2}\right)  ,
\]
and hence that $\mu^{2}+\nu^{2}\geq\frac{\delta}{3-\delta}\lambda^{2}$. Since
$\left\vert \mu\right\vert \geq\left\vert \nu\right\vert $ by assumption, we
have
\begin{align*}
P &  =\lambda^{2}\left(  \mu-\nu\right)  ^{2}+\mu^{2}\left(  \lambda
-\nu\right)  ^{2}+\nu^{2}\left(  \lambda-\mu\right)  ^{2}\\
&  \geq\lambda^{2}\left(  \mu-\nu\right)  ^{2}+\frac{\delta}{2\left(
3-\delta\right)  }\lambda^{2}\left(  \lambda-\nu\right)  ^{2}+\nu^{2}\left(
\lambda-\mu\right)  ^{2}.
\end{align*}
Now notice that
\[
\left\vert \operatorname*{Rm}\right\vert ^{2}=\lambda^{2}+\mu^{2}+\nu^{2}%
\leq3\lambda^{2}%
\]
and
\[
|\overset{\circ}{\operatorname*{Rm}}|^{2}=\frac{1}{3}\left[  \left(
\lambda-\mu\right)  ^{2}+\left(  \lambda-\nu\right)  ^{2}+\left(  \mu
-\nu\right)  ^{2}\right]  \leq\frac{4}{3}\left(  \lambda^{2}+\mu^{2}+\nu
^{2}\right)  \leq4\lambda^{2}.
\]
So if $\left\vert \nu\right\vert \leq\left\vert \lambda\right\vert /2$, we
have
\[
P\geq\frac{\delta}{2\left(  3-\delta\right)  }\lambda^{2}\left(  \lambda
-\nu\right)  ^{2}\geq\frac{\delta}{8\left(  3-\delta\right)  }\lambda^{4}%
\geq\frac{\delta}{96\left(  3-\delta\right)  }\left\vert \operatorname*{Rm}%
\right\vert ^{2}|\overset{\circ}{\operatorname*{Rm}}|^{2},
\]
while if $\left\vert \nu\right\vert >\left\vert \lambda\right\vert /2$, we
get
\begin{align*}
P &  >\lambda^{2}\left(  \mu-\nu\right)  ^{2}+\frac{\delta}{2\left(
3-\delta\right)  }\lambda^{2}\left(  \lambda-\nu\right)  ^{2}+\frac{1}%
{4}\lambda^{2}\left(  \lambda-\mu\right)  ^{2}\\
&  \geq\frac{\delta}{2\left(  3-\delta\right)  }\lambda^{2}\left[  \left(
\mu-\nu\right)  ^{2}+\left(  \lambda-\nu\right)  ^{2}+\left(  \lambda
-\mu\right)  ^{2}\right]  \geq\frac{\delta}{2\left(  3-\delta\right)
}\left\vert \operatorname*{Rm}\right\vert ^{2}|\overset{\circ}%
{\operatorname*{Rm}}|^{2},
\end{align*}
because $\frac{1}{4}>\frac{\delta}{2\left(  3-\delta\right)  }$.
\end{proof}


\begin{thebibliography}{9999999}                                                                                          %


\bibitem[ B-60]{B-60}M. Brown, \emph{A proof of the generalized Schoenflies
theorem, }Bull. Amer. Math. Soc. \textbf{66} (1960) 74-76.

\bibitem[CG-72]{CG-72}J. Cheeger and D. Gromoll, \emph{On the structure of
complete manifolds of nonnegative curvature,} Ann. of Math. (2) \textbf{96}
(1972) 413--443.

\bibitem[ C-44]{C-44}S.-S. Chern, \emph{A simple intrinsic proof of the
Gauss-Bonnet formula for closed Riemannian manifolds, }Ann. of Math.
\textbf{45} (1944) 747-752.

\bibitem[ CKE]{CKE}B. Chow, D. Knopf, et al., \emph{The Ricci flow on
3-manifolds} (e.g., \textit{Type I singularities }chapter,) book in progress.
http://www.math.uiowa.edu/\symbol{126}dknopf/

\bibitem[ D-94]{D-94}G. Drees, \emph{Asymptotically flat manifolds of
nonnegative curvature, }Diff. Geom. Appl. \textbf{4 }(1994) 77-90.

\bibitem[ ESS-89]{ESS-89}J. H. Eschenburg, V. Schroeder and M. Strake,
\emph{Curvature at infinity of open nonnegatively curved manifolds,} J. Diff.
Geom. \textbf{30} (1989) 155-166.

\bibitem[ G-97]{G-97}R.\ E.\ Greene, \textquotedblleft A genealogy of
noncompact manifolds of nonnegative curvature: history and
logic,\textquotedblright\ in \emph{Comparison Geometry, }ed.\ Grove and
Petersen, MSRI Publ. \textbf{30} (1997) 99--134.

\bibitem[GW-74]{GW-74}R.\ E.\ Greene and H. Wu, \emph{Integrals of subharmonic
functions on manifolds of nonnegative curvature,} Invent. Math. \textbf{27}
(1974) 265-298.

\bibitem[GW-76]{GW-76}R.\ E.\ Greene and H. Wu, $\emph{C}^{\infty}%
$\emph{-convex functions and manifolds of positive curvature,} Acta Math.
\textbf{137} (1976) 209-245.

\bibitem[ GW-82]{GW-82}R.\ E.\ Greene and H. Wu, \emph{Gap theorems for
noncompact Riemannian manifolds,} Duke. Math. J. \textbf{49} (1982) 731-756.

\bibitem[ GM-69]{GM-69}D.\ Gromoll and W.\ Meyer, \emph{On complete manifolds
of positive curvature, }Ann.\ of Math.\ \textbf{90} (1969) 75--90.

\bibitem[ H-93a]{H-93a}R. S. Hamilton, \emph{The Harnack estimate for the
Ricci flow,} J. Diff. Geom. \textbf{37 }(1993) 225-243.

\bibitem[ H-93b]{H-93b}R. S. Hamilton, \emph{Eternal solutions to the Ricci
flow,} J. Diff. Geom. \textbf{38 }(1993) 1-11.

\bibitem[ H-95a]{H-95a}R.\ S.\ Hamilton, \textquotedblleft The formation of
singularities in the Ricci flow,\textquotedblright\ in\emph{ Surveys in
Differential Geometry }\textbf{2} (1995) International Press, 7--136.

\bibitem[ H-95b]{H-95b}R.\ S.\ Hamilton, \emph{A compactness property for
solutions of the Ricci flow, }Amer.\ J.\ Math.\ \textbf{117 }(1995) 545--572.

\bibitem[ H-97]{H-97}R.\ S.\ Hamilton, \emph{Four-manifolds with positive
isotropic curvature, }Comm.\ Anal.\ Geom.\ \textbf{5} (1997) 1--92.

\bibitem[Hu-90]{Hu-90}G. Huisken, \emph{Asymptotic behavior for singularities
of the mean curvature flow,} J. Diff. Geom. \textbf{31} (1990) 285--299.

\bibitem[ I-93]{I-93}T.\ Ivey, \emph{Ricci solitons on compact }%
$3$\emph{-manifolds, }Differential Geom.\ Appl.\ \textbf{3} (1993) 301--307.

\bibitem[ K-88]{K-88}A. Kasue, \emph{A compactification of a manifold with
asymptotically nonnegative curvature, }Ann. Sci. Ec. Norm. Sup. IV Ser.
\textbf{21} (1988) 593-622.

\bibitem[KS-87]{KS-87}A. Kasue and K. Sugahara, \emph{Gap theorems for certain
submanifolds of Euclidean spaces and hyperbolic space forms, }Osaka J. Math.
\textbf{24} (1987) 679-704.

\bibitem[LT-87]{LT-87}P. Li and L.-F. Tam, \emph{Positive harmonic functions
on complete manifolds with nonnegative curvature outside a compact set,} Ann.
of Math. \textbf{125} (1987) 171--207.

\bibitem[M-59]{M-59}B. Mazur, \emph{On embeddings of spheres, }Bull. Amer.
Math. Soc. \textbf{65} (1959) 59-65.

\bibitem[PT-01]{PT-01}A. Petrunin and W. Tuschmann, \emph{Asymptotical
flatness and cone structure at infinity,} Math. Ann. \textbf{321} (2001) 775--788.

\bibitem[ Z-97]{Z-97}S. Zhu, \textquotedblleft The comparison geometry of
Ricci curvature,\textquotedblright\ in \emph{Comparison Geometry, }ed.\ Grove
and Petersen, MSRI Publ. \textbf{30} (1997) 221-262.
\end{thebibliography}
\end{document}